\def\no{\if01}
\def\iftwelvept{\no}

\def\ifusepdf{\no}
\def\ifpsfont{\no}

\iftwelvept
\documentclass[leqno,12pt]{amsart}
\else
\documentclass[leqno,10pt]{amsart}
\fi
\usepackage{amssymb}
\usepackage{amscd}
\usepackage{latexsym}
\usepackage{verbatim}
\usepackage[all]{xy}

\ifusepdf
\usepackage{hyperref}
\else\fi
\ifpsfont
\usepackage[T1]{fontenc}
\usepackage{times}
\else\fi

\iftwelvept
\else\fi

\theoremstyle{plain}
\newtheorem{Theorem}{Theorem}[section]

\newtheorem{Proposition}[Theorem]{Proposition}
\newtheorem{Lemma}[Theorem]{Lemma}
\newtheorem{Corollary}[Theorem]{Corollary}

\theoremstyle{definition}

\newtheorem{Definition}[Theorem]{Definition}
\newtheorem{Remark}[Theorem]{Remark}
\newtheorem{Construction}[Theorem]{Construction}

\newcommand{\Mfldf}{\mathsf{Mfld}^{\textup{fr}}}
\newcommand{\Diskf}{\mathsf{Disk}^{\textup{fr}}}

\newcommand{\ZZ}{\mathbf{Z}}

\newcommand{\RRR}{\mathbf{R}}
\newcommand{\CC}{\mathbb{C}}

\newcommand{\NNNN}{\operatorname{N}}
\newcommand{\GG}{\mathcal{G}}
\newcommand{\HH}{\operatorname{\mathcal{HH}}}

\newcommand{\ZZZ}{\mathbb{Z}}

\newcommand{\MM}{\mathbb{M}}

\newcommand{\uni}{\mathbf{1}}
\newcommand{\CCC}{\mathcal{C}}
\newcommand{\RPerf}{\operatorname{RPerf}}

\newcommand{\PR}{\operatorname{Pr}^{\textup{L}}}
\newcommand{\KKK}{\mathbf{K}}

\newcommand{\GGG}{\mathbb{G}}

\newcommand{\Hom}{\operatorname{Hom}}

\newcommand{\Ker}{\operatorname{Ker}}

\newcommand{\Comp}{\operatorname{Comp}}

\newcommand{\Perf}{\operatorname{Perf}}

\newcommand{\SP}{\operatorname{Sp}}
\newcommand{\SPS}{\operatorname{Sp}^{\Sigma}}

\newcommand{\Mod}{\operatorname{Mod}}

\newcommand{\SSS}{\mathcal{S}}

\newcommand{\colim}{\operatorname{colim}}
\newcommand{\CAT}{\operatorname{Cat}}
\newcommand{\Cat}{\textup{Cat}_{\infty}}

\newcommand{\Map}{\operatorname{Map}}

\newcommand{\Fun}{\operatorname{Fun}}
\newcommand{\Alg}{\operatorname{Alg}}

\newcommand{\End}{\operatorname{End}}

\newcommand{\Grp}{\operatorname{Grp}}

\newcommand{\CAlg}{\operatorname{CAlg}}

\newcommand{\LM}{\operatorname{\mathcal{LM}}}
\newcommand{\RM}{\operatorname{\mathcal{RM}}}

\newcommand{\Ind}{\operatorname{Ind}}

\newcommand{\Free}{\operatorname{Free}}

\newcommand{\assoc}{\operatorname{As}}

\newcommand{\eone}{\mathbf{E}_1}
\newcommand{\etwo}{\mathbf{E}_2}
\newcommand{\eenu}{\mathbf{E}_n}

\newcommand{\SO}{\textup{SO}}

\newcommand{\KS}{\mathbf{KS}}
\newcommand{\Lie}{\textup{Lie}}

\newcommand{\LMod}{\operatorname{LMod}}
\newcommand{\RMod}{\operatorname{RMod}}

\newcommand{\ST}{\operatorname{\mathcal{S}t}}

\newcommand{\kep}{k[\epsilon]}

\newcommand{\Proof}{{\sl Proof.}\quad}
\newcommand{\QED}{{\unskip\nobreak\hfil\penalty50\quad\null\nobreak\hfil
{$\Box$}\parfillskip0pt\finalhyphendemerits0\par\medskip}}

\begin{document}

\title{Bogomolov-Tian-Todorov theorem for Calabi-Yau categories}

\author{Isamu Iwanari}






\address{Mathematical Institute, Tohoku University, Sendai, Miyagi, 980-8578,
Japan}

\email{isamu.iwanari.a2@tohoku.ac.jp}

\maketitle

\section{Introduction}

The theorem of Bogomolov-Tian-Todorov \cite{Bo}, \cite{Tian}, \cite{Todorov}
says that
deformations of a compact Calabi-Yau manifold are unobstructed.
Namely,
the local moduli space is smooth. This unobstructedness is described as
a certain abelian property of the differential graded Lie algebra consisting of
derived global sections of the tangent sheaf.
In this note we consider a generalization to Calabi-Yau categories.
For a stable $\infty$-category $\CCC$ over a field of characteristic zero
(or a pretriangulated differential graded category), 
the deformation problem of $\CCC$ is controlled by
the shifted Hochschild cochain complex $\GG_{\CCC}:=\HH^\bullet(\CCC/k)[1]$, that
is a differential graded Lie algebra \cite[5.3.16]{DAG}
(see {\it loc. cit.} for the precise formulation in which it is formulated as a deformation problem on algebras over the little disks operad).
The purpose of this note is to prove a generalization of Bogomolov-Tian-Todorov theorem to Calabi-Yau categories (see Theorem~\ref{generalBTT}
and Corollary~\ref{CYBTT}):

\begin{Theorem}
let $k$ be a field of characteristic zero.
Let $\CCC$ be a $k$-linear proper Calabi-Yau stable $\infty$-category. 
Let $\GG_{\CCC}$ denote the shifted Hochschild cochain complex $\HH^\bullet(\CCC/k)[1]$
that is a differential graded (dg)
Lie algebra. Then $\GG_{\CCC}$ is quasi-isomorphic to
an abelian dg Lie algebra, that is, a dg Lie algebra with the zero bracket.
\end{Theorem}

We would like to take a brief look at a history
of work on Bogomolov-Tian-Todorov theorem related to the above theorem.
The above result was originally announced
in Katzarkov-Kontsevich-Pantev \cite[4.4.1]{KKP} for
the case of a proper Calabi-Yau
differential graded (dg) algebra $A$. Unfortunately, the proof has not appeared
at the writing of this note.
In the unpublished preprint \cite{IP} the author gave the proof of the above result for a proper Calabi-Yau
dg algebra $A$.
Subsequently, Tu \cite{Tu} proved a version of a
$\ZZ/2\ZZ$-graded finite-dimensional smooth cyclic dg ($A_\infty$) algebra under the assumption of the degeneration of
Hodge-de Rham spectral sequence in this situation.

\vspace{1mm}

Section~\ref{presection}
collects conventions and some of the notation that we will
use. We also review
some of background materials.
Sections~\ref{Periodmap}, \ref{contractionsection}, \ref{completion} are devoted to the proof of 
main results.
We may think that Section~\ref{Periodmap} is the first step,
Section~\ref{contractionsection} is the second step, and Section~\ref{completion}
consists of the completion of the proof.
We briefly mention some of the ingredients for the proof.
We will apply
a conceptual construction \cite{ID} of the algebraic
structure on the pair 
of the Hochschild cohomology spectrum (Hochschild cochain complex in the dg setting)
and the Hochschild homology spectrum (Hochschild chain complex)
$(\HH^\bullet(\CCC/k),\HH_\bullet(\CCC/k))$
of a $k$-linear stable $\infty$-category $\CCC$ (on a personal level, at the writing of \cite{IP} the author was dimly aware of ideas in the present arguments).
This algebraic structure is encoded by the Kontsevich-Soibelman operad.
Besides, we make use of the theory of factorization homology \cite{AF}
and the degeneration of Hodge-de Rham spectral sequence \cite{Kal0}, \cite{Ma}
for a smooth and proper stable $\infty$-category.
Roughly speaking, our approach is modeled on the work of Iacono
and Manetti \cite{IaMa}, who give an algebraic proof of the classical Bogomolov-Tian-Todorov theorem.
Along with the above inputs,
observations on
categorical phenomena, which do not appear in the classical situation, such as Lemmas~\ref{mysteriousisom} and~\ref{extensiontokernel},
Lemma~\ref{null}, and
Proposition~\ref{desiredcompatibility} play an important role.
We also hope that the reader finds our method natural
at least compared to that in \cite{IP}. For example, the argument of the present note does not require a deep result of the formality
of the Kontsevich-Soibelman operad and complicated computations
of $L_\infty$-algebras needed in  \cite{IP}.

\vspace{2mm}

{\it Acknowledgement.}
The author would like to thank A. Bondal, H. Minamoto,
T. Sano and A. Takahashi for discussions related to the subject of this note.
The author is partly supported by Grant-in-Aid for Young Scientists JSPS.

\section{Preliminaries}
\label{presection}
Unless otherwise stated, $k$ is a base field of characteristic zero.

\subsection{}

\label{Pre1}

{\it $(\infty,1)$-categories.}
We use 
the language of $(\infty,1)$-categories.
In several parts of this note one can work only with
homological algebras formulated in the theory
of model categories and their homotopy categories. But the author
feels that the theory of $(\infty,1)$-categories provides more consistent
and useful language for our purpose.
We adopt the theory of {\it quasi-categories} as a model of $(\infty,1)$-categories.
A quasi-category is a simplicial set which
satisfies the weak Kan condition of Boardman-Vogt.
Following \cite{HTT}, we shall refer to quasi-categories
as {\it $\infty$-categories}.
Our main references are \cite{HTT}
 and \cite{HA}.
 To an ordinary category, we can assign an $\infty$-category by taking
its nerve, and therefore
when we treat ordinary categories we often omit the nerve $\NNNN(-)$
and directly regard them as $\infty$-categories.
Here is a list of (some) of the conventions and notation that we will use:

\begin{itemize}

\item $\CCC^{op}$: the opposite $\infty$-category of an $\infty$-category. We also use the superscript ``op" to indicate the opposite category for ordinary categories and enriched categories.


\item $\SSS$: the $\infty$-category of $\infty$-groupoids/spaces.

\item $\SP$: the $\infty$-category of spectra.

\item $\Fun(A,B)$: the function complex for simplicial sets $A$ and $B$. If $A$ and $B$ are $\infty$-categories, we regard $\Fun(A,B)$ as the functor category.

\item $\Map_{\mathcal{C}}(C,C')$: the mapping space from an object $C\in\mathcal{C}$ to $C'\in \mathcal{C}$ where $\mathcal{C}$ is an $\infty$-category.
We usually view it as an object in $\mathcal{S}$ (cf. \cite[1.2.2]{HTT}).

\item 
$\Ind(\CCC)$: the $\infty$-category of Ind-objects (cf. \cite[5.3.5]{HTT}).

\item $S^1$: the circle. By default, we regard $S^1$ as an object of the $\infty$-category $\SSS$ of
spaces/$\infty$-groupoids.
When we think of $S^1$ as the group object, we denote it by $\SO(2)$
except we write $BS^1$ for $B\textup{SO}(2)$.

\end{itemize}

We fix the notation for the categories which we will frequently use.

\begin{itemize}

\item $\Mod_k$: Let $\Comp(k)$ denote the category
of chain complexes of $k$-vector spaces.
This category has a symmetric monoidal strucuture
given by the tensor product of chain complexes.
By inverting quasi-isomorphisms (see below), we have
the symmetric monoidal $\infty$-category $\Comp(k)[W^{-1}]$.
Regarding $k$ as a commutative ring spectrum,
we consider the symmetric monoidal $\infty$-category
$\Mod_k(\SP)$ of $k$-module objects in the $\infty$-category $\SP$ of 
spectra in the $\infty$-categorical formulation
(cf. \cite[3.3.3]{HA}).
Thanks to \cite[7.1.2.13]{HA}, there exists
a canonical equivalence
$\Mod_k(\SP)\simeq \Comp(k)[W^{-1}]$.
By abuse of notation, we write $\Mod_k$ for both $\infty$-categories.

\item $\Alg^{dg}(k)$: the combinatorial model category of associative (possibly noncommutative)
algebras such that the class of weak equivalences (resp. fibrations)
consists of quasi-isomorphisms (resp. termwise surjective maps)
when they are regarded as the underlying chain complex
(see e.g. \cite[7.1.4.5]{HA}).
Similarly, $\CAlg^{dg}(k)$ denotes the combinatorial model category of (genuine) commutative dg algebras over $k$
such that the class of weak equivalences (resp. fibrations)
consists of quasi-isomorphisms (resp. termwise surjective maps)
when they are regarded as underlying chain complexes
(cf. \cite[7.1.4.7, 7.1.4.10]{HA}).

\item $\Comp(A)$: Given $A\in \Alg^{dg}(k)$, we denote by
$\Comp(A)$ the combinatorial model category of dg left $A$-modules.
The class of weak equivalences consists of quasi-isomorphisms,
and the class of fibration consists of termwise surjective maps (see e.g. \cite[Theorem 3.3]{Six} or \cite[4.3.3.15]{HA}). 
We refer to it as the projective model structure.

\item $\Mod_A$: Given $A\in \CAlg^{dg}(k)$,
we write $\Mod_A$ for the $\infty$-category 
obtained from $\Comp(A)$ by inverting weak equivalences.
Let us regard $A$ as the image under $\CAlg^{dg}(k)=\CAlg(\Comp(k))\to \CAlg(\Mod_k)$
where the latter category denotes the $\infty$-category of commutative algebra objects
in $\Mod_k$.
If we write $\Mod_A(\Mod_k)$ for the $\infty$-category
of $A$-module objects in $\Mod_k$
in the $\infty$-categorical sense (see \cite[3.3.3]{HA}),
there exists a canonical equivalence
$\Mod_A\simeq \Mod_A(\Mod_k)$ (cf. \cite[4.3.3.17, 4.5.1.6]{HA}).

\item $\Alg(\Mod_k)$: the $\infty$-category obtained from $\Alg^{dg}(k)$
by inverting weak equivalences.
This category is also defined as the $\infty$-category $\Alg_{1}(\Mod_k)$ of
the associative algebra objects in the symmetric monoidal $\infty$-category
$\Mod_k$, i.e., algebras over the associative operad (see \cite[4.1.1]{HA}):
there exists a canonical equivalence $\Alg(\Mod_k)\simeq \Alg_1(\Mod_k)$
(see \cite[4.1.8.7]{HA}).
We refer to an object of $\Alg(\Mod_k)$ as a dg algebra over $k$.

\item $\LMod_A$: Given $A\in \Alg^{dg}(k)$,
$\LMod_A$ is the $\infty$-category obtained from $\Comp(A)$
by inverting weak equivalences. We refer to this category as the $\infty$-category
of left dg $A$-modules.
Similarly, the $\infty$-category $\RMod_A$ of right dg $A$-modules
is defined to be the $\infty$-category obtained from $\Comp(A^{op})$.
We let $A$ denote the image of $A$ in $\Alg_1(\Mod_k)$ (by abuse of notation)
and let $\LMod_A(\Mod_k)$ denote the $\infty$-category of left $A$-module objects in $\Mod_k$ (cf. \cite[4.1.1.1]{HA}).
According to  \cite[4.3.3.17]{HA}, there exists
 a canonical equivalence $\LMod_A\simeq \LMod_{A}(\Mod_k)$.
 If we define $\RMod_A(\Mod_k)$ in a similar way,
we have $\RMod_A\simeq \RMod_{A}(\Mod_k)$.

\item $\Lie_k$: the $\infty$-category of dg Lie algebras over $k$.
We use the Lie operad
defined as a $\Mod_k$-enriched $\infty$-operad developed
in \cite{Hau1} and \cite{Hau2}.
For a $k$-linear presentable symmetric monoidal stable $\infty$-category,
i.e., an object of $\CAlg(\PR_k)$ (see Section~\ref{Pr3})
(or a small $k$-linear symmetric monoidal stable $\infty$-category) $\mathcal{M}^\otimes$,
we denote by $\Lie(\mathcal{M}^\otimes)$ (or $\Lie(\mathcal{M})$) 
the $\infty$-category of algebras in $\mathcal{M}^\otimes$
over the Lie operad (cf. \cite[Section 3]{Hau2}).

Let $\textup{Lie}^{dg}(k)$ be the usual category of (genuine) dg Lie algebras over $k$. This category admits a combinatorial
model structure such that the class of weak equivalences consists of
quasi-isomorphisms (when regarded as underlying complexes),
(see \cite[Theorem 4.1.1]{Hin}).
We define $\Lie_k$ to be the $\infty$-category obtained from $\textup{Lie}^{dg}(k)$ by
 inverting weak equivalences.
 By the comparison  \cite[Theorem 4.10, Corollary 4.11]{Hau2}
 between algebras over enriched $\infty$-operads and algebras over the classical dg operads
and the comparison \cite[Theorem 4.7.2]{Hin} between algebras over quasi-isomorphic
dg operads in characteristic zero, there exists a canonically defined
equivalence of $\infty$-catrgories $\Lie_k\simeq \Lie(\Mod_k)$.

\end{itemize}

{\it From model categories to $\infty$-categories.}
For the reader's convenience, we recall Lurie's construction
by which one can obtain
$\infty$-categories
from a category (more generally
$\infty$-category) endowed with a prescribed collection of
morphisms (see \cite[1.3.4, 4.1.7, 4.1.8]{HA} for details).
It can be viewed as an alternative approach to Dwyer-Kan
hammock localization.
Let $\mathcal{D}$ be a category and let $W$ be a collection
of morphisms in $\mathcal{D}$ which is closed under composition
and contains all isomorphisms.
For $(\mathcal{D},W)$, there is an $\infty$-category $\NNNN(\mathcal{D})[W^{-1}]$
and a functor
$\xi:\textup{N}(\mathcal{D})\to \NNNN(\mathcal{D})[W^{-1}]$
such that for any $\infty$-category $\CCC$ the composition
induces a fully faithful functor $\Fun(\NNNN(\mathcal{D})[W^{-1}],\CCC)\to \Fun(\textup{N}(\mathcal{D}),\CCC)$
whose essential image consists of those functors $F:\textup{N}(\mathcal{D})\to\CCC$ such that $F$ carry morphisms lying in $W$
to equivalences in $\CCC$.
We shall refer to $\NNNN(\mathcal{D})[W^{-1}]$
as the $\infty$-category obtained from $\mathcal{D}$
by inverting morphisms in $W$.

\subsection{$\mathbf{E}_n$-algebras and universal enveloping algebras.}
\label{Pre2}

We use the theory of algebras over little cubes (disks) operads.
Among various formulations, 
we adopt the $\infty$-operad of little $n$-cubes, which we denote
by $\mathbf{E}^\otimes_n$ (cf. \cite[5.1]{HA}).
This choice is not essential.
For a symmetric monoidal $\infty$-category $\CCC^\otimes$,
we write $\Alg_{n}(\CCC)$ for
the $\infty$-category of algebra objects over $\mathbf{E}^\otimes_n$ in $\CCC^\otimes$.
We refer to an object of $\Alg_{n}(\CCC)$ as an $\eenu$-algebra in $\CCC$.
If we denote by $\assoc^\otimes$ the associative operad (\cite[4.1.1]{HA}),
there is the standard equivalence $\assoc^\otimes\simeq \eone^\otimes$
of $\infty$-operads. We usually identify $\Alg_{1}(\CCC)$
with the
$\infty$-category of associative algebras in $\CCC$ (this notation
is consistent with $\Alg_1(\Mod_k)$).
In this note, we need only $\Alg_1(\CCC)$ and $\Alg_2(\CCC)$.

\vspace{2mm}

Next we review the universal enveloping algebra functor $U_n:\Lie_k\to \Alg_n(\Mod_k)$.
Let $\Omega_L:\Lie_k\to \Grp(\Lie_k)$ be the based loop space functor
which carries $L$ to the fiber product $0\times_L0$ in $\Lie_k$
where $0$ is the ``zero'' dg Lie algebra (the construction of $\Omega_L$
is left to the reader). Here $\Grp(\Lie_k)$ is the $\infty$-category
of group objects (i.e., group-like $\eone$-monoids) (see \cite[7.2.2.1]{HTT}).
Consider the adjoint pair
\[
\xymatrix{
B_L:\Grp(\Lie_k)  \ar@<0.5ex>[r] &  \Lie_k:\Omega_L  \ar@<0.5ex>[l]  
}
\]
where $B_L$ is the left adjoint which carries a group object $G$
to the geometric realization of the simplicial object determined by $G$
in $\Lie_k$.
As explained in \cite[Chap. 6, 1.6.4]{Gai2},
$B_L$ and $\Omega_L$ are mutually inverse equivalences.
Moreover, if we forget the Lie algebra structures (and group structures),
$\Omega_L$ (resp. $B_L$) gives rise to the shift functor $[-1]$ (resp. $[1]$)
at the level of underlying
complexes.
By \cite[Chap. 6, 1.6.2]{Gai2},
the inclusion $\Grp(\Lie_k)\subset \textup{Mon}(\Lie_k)$ into the $\infty$-category
of monoid objects is an equivalence.

Consider the standard adjoint pair
$U_1:\Lie_k\rightleftarrows \Alg_1(\Mod_k):res_{\eone/Lie}$,
where $U_1$ is the ``standard'' universal enveloping algebra functor.
As proved in \cite[Chap. 6, 6.1.2]{Gai2},
$U_1$ is equivalent to
the composite $\Lie_k\stackrel{\Omega_L}{\longrightarrow} \Grp(\Lie_k)\simeq \Alg_1(\Lie_k^\times)\stackrel{\Alg_1(Ch_\bullet)}{\longrightarrow}\Alg_1(\Mod_k)$.
Here $\Alg_1(\Lie_k^\times)$ be the $\infty$-category of associative algebras
in the Cartesian monoidal $\infty$-category $\Lie_k^\times$
(that is, the $\infty$-category of the monoid objects),
and $Ch_\bullet:\Lie_k\to \Mod_k$ is the Chevalley-Eilenberg chain functor.

\begin{Lemma}
\label{identifyI2}
Let $\Grp^{n-1}(\Lie_k^\times)=\Grp(\Grp(\cdots \Grp(\Lie_k)\cdots ))$
be the $\infty$-category of $(n-1)$-fold group objects
and let $\Alg_{n-1}(\Lie_k^\times)$ be the $\infty$-category
of $\mathbf{E}_{n-1}$-algebra objects in $\Lie_k$
endowed with the symmetric monoidal structure given by the
Cartesian product. 
There exists a sequence of adjoint pairs
\[
\xymatrix{
\Mod_k \ar@<0.5ex>[r]^{\textup{Free}_{\Lie}} & \Lie_k \ar@<0.5ex>[r]^(0.3){\Omega^{n-1}_L} \ar@<0.5ex>[l]^{res_{Lie/M}} & \Grp^{n-1}(\Lie_k^\times)\simeq \Alg_{n-1}(\Lie_k^{\times}) \ar@<0.5ex>[r] \ar@<0.5ex>[l]^(0.7){B^{n-1}_L} & \Alg_{n-1}(\Alg_1(\Mod_k))\simeq \Alg_n(\Mod_k) \ar@<0.5ex>[l]
}
\]
such that
\begin{enumerate}
\item $res_{Lie/M}$ is the forgetful functor, and
 $\textup{Free}_{Lie}$ is the free Lie algebra functor,

\item $\Omega_L^{n-1}$ and $B_L^{n-1}$ are mutually inverse equivalences,

\item the right adjoint pair is given by
$(\Alg_{n-1}(U_1),\Alg_{n-1}(res_{\eone/Lie}))$
induced by the adjoint pair $(U_1,res_{\eone/Lie})$,

\item $res_{Lie/M}\circ B_L^{n-1}$ is equivalent to the compostion
of the forgetful functor $\Grp^{n-1}(\Lie_k)\to \Mod_k$ and 
the shift functor $[n-1]:\Mod_k\to \Mod_k$.
In particular, the composite $\Alg_n(\Mod_k)\to \Mod_k$ of right adjoint funtors
is equivalent to the forgetful functor $\Alg_n(\Mod_k)\to \Mod_k$ followed by
the shift functor $[n-1]$.
\end{enumerate}
\end{Lemma}

\begin{Remark}
By ease of notation, we often write $U_1$ for $\Alg_{n-1}(U_1)$. 
\end{Remark}

\Proof
As reviewed above, the result
\cite[Chap. 6,  Prop.1.6.4]{Gai2}
says that $B_L$ and $\Omega_L$ are mutually inverse equivalences.
Thus (2) follows.
Moreover, it follows from the proof of {\it loc. cit.}
 that $res_{Lie/M}\circ B_L^{n-1}$ is equivalent to the composition
of the forgetful functor $\Grp^{n-1}(\Lie_k)\to \Mod_k$ and 
the shift functor $[n-1]:\Mod_k\to \Mod_k$.
The middle equivalence 
follows from \cite[Chap. 6, Lemma 1.6.2]{Gai2}.
The right equivalence is the Dunn additivity theorem \cite[5.1.2.21]{HA}.
Other adjoint pairs are standard ones.
\QED

\begin{Definition}
The universal enveloping algebra functor $U_n:\Lie_k\to \Alg_n(\Mod_k)$
is defined to be the left adjoint to $B_L^{n-1}\circ \Alg_{n-1}(res_{\eone/Lie})$.
That is, $U_n=\Alg_{n-1}(U_1)\circ \Omega_{L}^{n-1}$ (up to the right equivalence
given by the Dunn additivity).
\end{Definition}

\begin{Remark}
This definition of $U_n$ is different from the definition we give in \cite[Section 3.5]{IM}.
But it is not difficult to prove that both definitions of $U_n$ coincide (use \cite[Propoition 3.3]{IM}).
In this note, we will not use this comparison so that we omit the proof.
\end{Remark}

\begin{Remark}
In this note we will use only $U_1$ and $U_2$.
\end{Remark}

\subsection{Linear categories}
\label{Pr3}
Suppose that $A\in \CAlg^{dg}(k)$.
We will fix our convention on $A$-linear stable $\infty$-categories.
More detailed account can be found in \cite[Section 3.1]{IM}.
Let $\PR_{\textup{St}}$ be the $\infty$-category
of stable presentable $\infty$-categories whose morphisms are 
functors which preserve small colimits (see \cite[5.5.3.1]{HTT}, \cite[4.8.2]{HA}).
It admits a symmetric monoidal structure defined in \cite[4.8.1.15]{HA}.
We can think of the symmetric monoidal $\infty$-category $\Mod_A^\otimes$
as a commutative algebra object in $\PR_{\textup{St}}$ with respect to this
symmetric monoidal strucuture. We define the $\infty$-category $\PR_A$
of $A$-linear (stable) presentable $\infty$-categories to be
the $\infty$-category $\Mod_{\Mod_A^\otimes}(\PR_{\textup{St}})$
of $\Mod_A^\otimes$-module objects in $\PR_{\textup{St}}$.
Let $\ST$ be the $\infty$-category of small stable idempotent-complete
$\infty$-categories whose functors are spanned by exact functors
(see e.g. \cite[Section 2.2]{BGT1}).
There is a symmetric monoidal strucutre on $\ST$ (see e.g. \cite[Section 5]{ID},
\cite[Section 3.1]{IM}).
If we write $\Perf_A^\otimes$ for the stable symmetric monoidal subcategory of
$\Mod_A^\otimes$ spanned by compact objects (equivalently, dualizable objects),
then $\Perf_A^\otimes$ can be regarded as a commutative algebra object in $\ST$.
The $\infty$-category $\ST_A$ of small $A$-linear 
stable $\infty$-categories is defined to be
the $\infty$-category of $\Perf_A^\otimes$-module objects in $\ST$.
The $\infty$-categories $\ST$ and $\ST_A$ can also
be defined in terms of enriched categories
such as dg categories and spectral categories and their Morita equivalences
(see \cite{Co}).

We recall the notion of properness and smoothness.
Suppose that $\CCC\in \ST_k$. We say that
$\CCC$ is proper over $k$ if the hom/mapping complex $\textup{Hom}_k(C,C')$, defined an object of $\Mod_k$,
is bounded with finite dimensional cohomology for any two objects $C,C'\in \CCC$.
We say that $\CCC$ is smooth over $k$
if the identity functor is a compact object in  $\mathcal{M}or_k(\Ind(\CCC),\Ind(\CCC))$ where $\mathcal{M}or_k(\Ind(\CCC),\Ind(\CCC))$
is the internal hom/mapping object in $\PR_k$ (see e.g. \cite[Lemma 5.1]{ID}).

\subsection{Hochschild homology and Hochschild cohomology}
\label{hochschildhomologysection}
We briefly recall Hochschild homology. We refer to \cite[Section 6]{ID} for details.

The Hochschild homology (chain) complex functor is defined as a symmetric monoidal functor 
\[
\HH_\bullet(-/k):\ST_k \longrightarrow \Mod_k^{\SO(2)}=\Fun(BS^1,\Mod_k)
\]
which carries $\mathcal{C} \in \ST_k$ to the Hochschild homology complex (Hochchild chain complex) $\HH_\bullet(\mathcal{C}/k)$.
Here the symmetric monoidal structure on $\Fun(BS^1,\Mod_k)$
is induced by that of $\Mod_k$ (cf. \cite[2.1.3.4]{HA}).
We refer the reader to \cite[Section 6, 6.14]{ID} for the construction
of $\HH_\bullet(-/k)$
(in {\it loc. cit.} we use the symbol $\HH_\bullet(\mathcal{C})$ instead of  $\HH_\bullet(\mathcal{C}/k)$).
Let $B\in \Alg_1(\Mod_k)$.
Let $\textup{LPerf}_B$ be the smallest stable idempotent-complete subcategory of $\LMod_B$
which contains $B$.
We will write  $\HH_\bullet(B/k)$ for $\HH_\bullet(\textup{LPerf}_B/k)$.

Let $\Theta_k:\Alg_1(\Mod_k)\to \PR_k$
the symmetric monoidal functor informally given by $B\mapsto \LMod_B$ (see e.g. \cite[4.8.3]{HA} for the construction or use the formalism of spectral categories).
Since $\Theta_k$ is a symmetric monoidal functor,
it gives rise to $\theta_k:\Alg_2(\Mod_k)\simeq \Alg_1(\Alg_1(\Mod_k))\stackrel{\Alg_1(\Theta_k)}{\longrightarrow} \Alg_1(\PR_k)$.
If $B$ is an $\etwo$-algebra in $\Mod_k$, then $\LMod_B$ has an associative monoidal
structure defined as $\theta_k(B)=\LMod_B^\otimes$.  
Moreover, $\textup{LPerf}_B$ is promoted to an object of $\Alg_1(\ST_k)$.
Thus, the symmetric monoidal functor $\HH_\bullet(-/k)$ gives $\HH_\bullet(B/k)$
which belongs to $\Alg_1(\Mod^{\SO(2)}_k)$.

We review the Hochschild cohomology of $\CCC\in \ST_k$.
Let $\Ind(\CCC)$ be the Ind-category that belongs to
$\PR_k$.
Let  $\mathcal{E}nd_k(\Ind(\CCC))\in \Alg_1(\PR_k)$
be the endomorphism algebra object of $\Ind(\CCC)$
which is endowed with a tautological action on $\Ind(\CCC)$ (see \cite[4.7.1]{HA} for the formulation).
By \cite[4.8.5.11, 4.8.5.16]{HA} there exists an adjoint pair
\[
\xymatrix{
\theta_k:\Alg_{2}(\Mod_k)  \ar@<0.5ex>[r] &    \Alg_{1}(\PR_k):E_k  \ar@<0.5ex>[l]  
}
\]
such that $E_k$ sends $\mathcal{M}^\otimes$ to the endomorphism algebra of the unit object $\uni_{\mathcal{M}}$ where we think of an object of
$\Alg_{1}(\PR_k)$ as an associative monoidal $k$-linear $\infty$-category.
The Hochschild cohomology complex (Hochschild cochain complex) $\HH^\bullet(\CCC/k)$ 
is defined to be $E_k(\mathcal{E}nd_k(\Ind(\CCC)))$.
We refer to  $\HH^\bullet(\CCC/k)$ as the Hochschild cohomology complex of $\CCC$ (over $k$).
Note that 
the counit map $\theta_k(E_k(\mathcal{E}nd_k(\Ind(\CCC))))\to \mathcal{E}nd_k(\Ind(\CCC))$ and the tautological action of 
$\mathcal{E}nd_k(\Ind(\CCC))$ on $\Ind(\CCC)$ determine
the left module action of $\LMod_{\HH^\bullet(\CCC/k)}^\otimes$ on $\Ind(\CCC)$.

Though we use the word ``homology'' and ``cohomology'',
$\HH_\bullet(\CCC/k)$ and $\HH^\bullet(\CCC/k)$
are not graded modules obtained by passing to (co)homology but spectra (or chain complexes) with algebraic structures.
Passing to homology of $\HH_\bullet(\CCC/k)$ and $\HH^\bullet(\CCC/k)$,
we have graded vector spaces which we denote by
$HH_*(\CCC/k)$ and $HH^*(\CCC/k)$, respectively.
We refer to $HH_*(\CCC/k)$ and $HH^*(\CCC/k)$
as Hochschild homology and Hochshchild cohomology, respectively.

\section{Hochschild pair and the classifying map}

\label{Periodmap}

The main result of this section is Proposition~\ref{mainreduction3}
which reduces the main theorem to the injectivity of
a certain map of graded vector spaces.
To construct the certain map (a sort of classifying map) we use the pair
$(\HH^\bullet(\CCC/k),\HH_\bullet(\CCC/k))$.

\subsection{}
\label{KSalgebrasection}
We briefly review the algebra of $(\HH^\bullet(\CCC/k),\HH_\bullet(\CCC/k))$
encoded by Kontsevich-Soibelman operad $\KS$.
For details, we refer the reader to \cite{ID}
where we construct an algebra $(\HH^\bullet(\CCC/k),\HH_\bullet(\CCC/k))$
over $\KS$. 
We do not recall the operad $\KS$. Instead,
we will give an equivalent formulation which is sutaible for our purpose. 
According to \cite[Theorem 1.2]{ID},
a $\KS$-algebra $(\HH^\bullet(\CCC/k),\HH_\bullet(\CCC/k))$ in $\Mod_k$ is equivalent to
giving the following triple:
an $\etwo$-algebra $\HH^\bullet(\CCC/k)\in \Alg_{2}(\Mod_k)$,
an $k$-module with $\SO(2)$-action $\HH_\bullet(\CCC/k)\in \Mod_k^{\SO(2)}$,
and a left $\HH_\bullet(\HH^\bullet(\CCC/k)/k)$-module $\HH_\bullet(\CCC/k)$
in $\Mod_k^{\SO(2)}$ that is an object of $\LMod_{\HH_\bullet(\HH^\bullet(\CCC/k)/k)}(\Mod_k^{\SO(2)})$
lying over $\HH_\bullet(\CCC/k)\in \Mod_k^{\SO(2)}$.
We sketch the construction of a left $\HH_\bullet(\HH^\bullet(\CCC/k)/k)$-module $\HH_\bullet(\CCC/k)$.
Recall the adjoint pair $(\theta_k,E_k)$ between $\Alg_2(\Mod_k)$ and $\Alg_1(\PR_k)$
from Section~\ref{hochschildhomologysection}.
Consider the associated counit map
\[
\LMod_{\HH^\bullet(\CCC/k)}=\theta_k\circ E_k(\mathcal{E}nd_k(\Ind(\CCC)))\to \mathcal{E}nd_k(\Ind(\CCC))
\]
in $\Alg_1(\PR_k)$
and the tautological left action of $\mathcal{E}nd_k(\Ind(\CCC))$ on $\Ind(\CCC)$.
We obtain a left $\LMod_{\HH^\bullet(\CCC/k)}$-module $\Ind(\CCC)$
in $\PR_k$.
Passing to full subcategories of compact objects we have
a left $\Perf_{\HH^\bullet(\CCC/k)}$-module $\CCC$ in $\ST_k$.
Applying $\HH_\bullet(-/k)$ to the left $\Perf_{\HH^\bullet(\CCC/k)}$-module $\CCC$
we get a left $\HH_\bullet(\HH^\bullet(\CCC/k)/k)$-module $\HH_\bullet(\CCC/k)$.
Let $\End(\HH_\bullet(\CCC/k))$ denote the endomorphism algebra object
of $\HH_\bullet(\CCC/k)$ in $\Mod_k^{\SO(2)}$, which is an object of $\Alg_{1}(\Mod_k^{\SO(2)})$.
The left module structure amounts to a morphism
\[
A_{\CCC}:\HH_\bullet(\HH^\bullet(\CCC/k)/k)\longrightarrow \End(\HH_\bullet(\CCC/k))
\]
in $\Alg_{1}(\Mod_k^{\SO(2)})$.

\begin{Remark}
\label{restrictedaction}
Note that $\LMod_{\HH^\bullet(\CCC/k)}$ is the smallest associative monoidal
subcategory of $\mathcal{E}nd_k(\Ind(\CCC))$
which is closed under small colimits and consists of the unit object (the identity functor).
Let $\CCC'$ be a stable idempotent-complete subcategory of $\CCC$.
Then by the restriction to $\Ind(\CCC')$
the action of the associative monoidal $\infty$-category
$\LMod_{\HH^\bullet(\CCC/k)}$ on $\Ind(\CCC)$ induces
an action on $\Ind(\CCC')$.
Consequently, it gives rise to the left action of $\HH_\bullet(\HH^\bullet(\CCC/k)/k)$
on $\HH_\bullet(\CCC'/k)$.
Moreover, the morphism
$\HH^\bullet(\CCC/k)\to \HH_\bullet(\HH^\bullet(\CCC/k)/k)$
of $\eone$-algebras induces the action of $\HH^\bullet(\CCC/k)$ on $\HH_\bullet(\CCC'/k)$.
\end{Remark}

\subsection{}
Let $\GG_{\CCC}$ be the shifted Hochschild cohomology/cochain complex,
that is defined to be the image of the $\etwo$-algebra $\HH^\bullet(\CCC/k)$
under the forgetful functor $res_{\etwo/Lie}=B_L\circ \Alg_1(res_{\eone/Lie}):\Alg_{2}(\Mod_k)\to \Lie_k$ (see Section~\ref{Pre2}).
The underlying complex of $\GG_{\CCC}$ is quasi-isomorphic/equivalent to 
$\HH^\bullet(\CCC/k)[1]$ (see Lemma~\ref{identifyI2} (4)).
The left $\HH_\bullet(\HH^\bullet(\CCC/k)/k)$-module $\HH_\bullet(\CCC/k)$
induces
actions of $\GG_{\CCC}$ and $\GG_{\CCC}^{S^1}$ on $\HH_\bullet(\CCC/k)$.
Here $\GG_{\CCC}^{S^1}$ means the cotensor with $S^1$ in $\Lie_k$
so that $\GG_{\CCC}^{S^1}\simeq \GG_{\CCC}\times_{\GG_{\CCC}\times \GG_{\CCC}}\GG_{\CCC}$  in $\Lie_k$.
The following is the strategy of the construction.
The counit map of the adjoint pair $(U_2,res_{\etwo/Lie})$
induces $U_2(\GG_{\CCC})=U_2(res_{\etwo/Lie}(\HH^\bullet(\CCC/k))\to \HH^\bullet(\CCC/k)$.
By applying $\HH_\bullet(-/k)$, it gives rise to
$\HH_\bullet(U_2(\GG_{\CCC})/k)\to \HH_\bullet(\HH^\bullet(\CCC/k)/k)$.
Consider the sequence in $\Alg_1(\Mod_k^{\SO(2)})$:
\[
U_1(\GG_{\CCC})\to U_1(\GG_{\CCC}^{S^1})\simeq \HH_\bullet(U_2(\GG_{\CCC})/k)\to \HH_\bullet(\HH^\bullet(\CCC/k)/k) \stackrel{A_{\CCC}}{\longrightarrow} \End(\HH_\bullet(\CCC/k))
\]
where the first morphism is determined by the cotensor with $S^1\to \ast$,
the second morphism is the equivalence which comes from the nonabelian
Poincar\'e duality of factorization homology (see the review before the proof of Lemma~\ref{middlesquare}). 
Since the $\SO(2)$-action on $U_1(\GG_{\CCC})$ is trivial, the composite of the sequence
amounts to
a morphism
$U_1(\GG_{\CCC})\longrightarrow \End(\HH_\bullet(\CCC/k))^{\SO(2)}$
in $\Alg_1(\Mod_k)$. By the adjoint pair $(U_1,res_{\eone/Lie})$, this morphism
gives us
\[
A_{\CCC}^L:\GG_{\CCC}\longrightarrow \End^L(\HH_\bullet(\CCC/k))^{\SO(2)}
\]
in $\Lie_k$.
Here $\End(\HH_\bullet(\CCC/k))$ is the endomorphism
algebra object which is defined as an object of $\Alg_1(\Mod_k^{\SO(2)})$,
and $\End^L(\HH_\bullet(\CCC/k))\in \Lie_k^{\SO(2)}=\Fun(BS^1,\Lie_k)$ is the dg Lie algebra obtained from $\End(\HH_\bullet(\CCC/k))\in \Alg_{1}(\Mod_k^{\SO(2)})$.
Let $\End^L(\HH_\bullet(\CCC/k))^{\SO(2)}\in \Lie_k$ be the homotopy fixed points of the $\SO(2)$-action.
Similarly, $U_1(\GG_{\CCC}^{S^1})\to  \End(\HH_\bullet(\CCC/k))$ gives rise to
\[
\widehat{A}_{\CCC}^{L}:\GG^{S^1}_{\CCC}\longrightarrow \End^L(\HH_\bullet(\CCC/k))
\]
in $\Lie_k^{\SO(2)}=\Fun(BS^1,\Lie_k)$.
Namely, we have the diagram
\[
\label{diagramA}
\vcenter{
\xymatrix{
\GG_{\CCC} \ar[r]^(0.3){A^L_{\CCC}} \ar[d]_{\textup{diagonal}} & \End^L(\HH_\bullet(\CCC/k)) \\
\GG_{\CCC}^{S^1} \ar[ur]_{\widehat{A}^L_{\CCC}} & 
}\tag{3.1}}
\]
in $\Fun(BS^1, \Lie_k)$, which commutes up to canonical homotopy.



\subsection{}

{\it Mixed complexes.}
We briefly review mixed complexes
from the perspective of models of chain complexes with $\SO(2)$-actions. 
Let $k[\epsilon]$ be the free commutative dg algebra
over $k$ which is generated by one element $\epsilon$ of homological
degree $1$. Namely, $k[\epsilon]$ is a graded algebra
$k\oplus k\cdot \epsilon$ endowed with the zero differential.
Consider the category $\Comp(k[\epsilon])$ of dg $\kep$-modules (see Section~\ref{Pre1}).
A dg $\kep$-module is often referred to as a mixed complex.
By the standard fact (cf. \cite{Hoy}), there exists
a canonical equivalence
\[
\Fun(BS^1, \Mod_k)\simeq \Mod_{\kep} =\Comp(\kep)[W^{-1}]
\]
of $\infty$-categories.

We consider the adjoint pair
\[
F:\Comp(\kep) \rightleftarrows \Comp(k):R
\]
where $F$ is the left adjoint given by the forgetful functor,
and the right adjoint $R$ is given by the tensor with $k[\eta]$.
Here $k[\eta]$ is the free commutative dg algebra generated by one element
$\eta$ of the homological degree $-1$ (i.e., 
$k[\eta]\simeq H^*(S^1,k)$)
while $k[\eta]$ is regarded as $\kep [-1]$ as a $\kep$-module
in this context.
Both $F$ and $R$ preserve quasi-isomorphisms.
By inverting quasi-isomorphims
we have the induced adjoint pair
$\Comp(\kep)[W^{-1}]\simeq \Fun(BS^1, \Mod_k)\rightleftarrows \Mod_k\simeq \Comp(k)[W^{-1}]$. The induced left adjoint is the forgetful functor, and
the induced right adjoint is given by $M\mapsto M^{S^1}$
where $M^{S^1}$ denotes the cotensor with $S^1$.

\vspace{2mm}

$\bullet$ In the rest of this section  we assume that
the $\SO(2)$-action on $\HH_\bullet(\CCC/k)$ is trivial.

\begin{Construction}
\label{modelconstruction}
Consider the diagram
 obtained from the diagram~(\ref{diagramA}) by forgetting Lie algebra structures.
We will construct its model in $\Comp(\kep)$. 

Suppose that $\GGG_{\CCC}$ is
an object of $\Comp(k)$ which represents
the image of $\GG_{\CCC}$ under $\Lie_k \to \Mod_k$.
Consider $R(\GGG_{\CCC})=\GGG_{\CCC}\otimes_kk[\eta]$,
which is a model of $\GG_{\CCC}^{S^1}$ in $\Comp(\kep)$.
Since $k[\eta]\simeq \kep[-1]$,
it follows that $R(\GGG_{\CCC})=\GGG_{\CCC}\otimes_kk[\eta]$
is a cofibrant object with respect to the projective model structure.
The canonical inclusion $\GGG_{\CCC}\to \GGG_{\CCC}\otimes_kk[\eta]$
induced by $k\to k[\eta]$ is a model of $\GG_{\CCC}\to \GG_{\CCC}^{S^1}$.

Consider 
$HH_*(\CCC/k)$ to be a complex with the zero differential.
Let $\End(HH_*(\CCC/k))$ be the endomorphism dg algebra
(that is indeed a graded algebra).
This dg algebra $\End(HH_*(\CCC/k))\in \Alg^{dg}(k)$ is a model of $\End(\HH_\bullet(\CCC/k))\in \Alg(\Mod_k)$
(Note that $\End(HH_*(\CCC/k))$ with the tautological map $\End(HH_*(\CCC/k))\otimes_kHH_*(\CCC/k)\to HH_*(\CCC/k)$ constitutes an endomorphism object of $\HH_\bullet(\CCC/k)$
in $\Mod_k$
because we work over a field $k$.
Thus it follows from \cite[4.7.1.40, 3.2.2.5]{HA} that the genuine dg algebra $\End(HH_*(\CCC/k))$ represents an endomorphism algebra of $\HH_\bullet(\CCC/k)$).

Consider $HH_*(\CCC/k)$ to be the $\kep$-module
induced by the restriction of the $k$-module structure $HH_*(\CCC/k)$ along
$\kep\to k$.
We refer to the $\kep$-module structure arising from a $k$-module by
restriction along $\kep\to k$ as the trivial $\kep$-module structure.
Consider 
$\End(HH_*(\CCC/k))$
to be the graded algebra with the trivial $\kep$-module structure.
This is a model of $\End(\HH_\bullet(\CCC/k))$.
Let $\End^L(HH_*(\CCC/k))$
be the associated dg Lie algebra (with the zero differential and the trivial $\kep$-module structure).
But we write $\End(HH_*(\CCC/k))$ for the underlying complex of 
$\End^L(HH_*(\CCC/k))$, obtained by forgetting the Lie algebra structure.
Since $\GGG_{\CCC}\otimes_kk[\eta]$ is cofibrant in $\Comp(\kep)$,
we can take the commutative diagram
\[
\label{diagramB}
\vcenter{
\xymatrix{
\GGG_{\CCC} \ar[r]^(0.3){\mathbb{A}^L_{\CCC}} \ar[d] & \End(HH_\ast(\CCC/k)) \\
\GGG_{\CCC}\otimes_kk[\eta] \ar[ur]_{\widehat{\mathbb{A}}_{\CCC}^L} & 
}\tag{3.2}}
\]
in $\Comp(\kep)$ which represents the diagram in $\Fun(BS^1,\Mod_k)$
 obtained from the diagram~(\ref{diagramA}) by forgetting Lie algebra structures.
Here we regard $ \End(HH_\ast(\CCC/k))$ as the underlying complex of
the dg Lie algebra $\End^L(HH_*(\CCC/k))$.
We define $\mathbb{I}_{\CCC}:\GGG_{\CCC}\otimes \eta\to \End(HH_\ast(\CCC/k))$
to be the composite $\mathbb{I}_{\CCC}:\GGG_{\CCC}\otimes \eta\hookrightarrow \GGG_{\CCC}\otimes_kk[\eta] \stackrel{\widehat{\mathbb{A}}_{\CCC}^L}{\longrightarrow} \End(HH_\ast(\CCC/k))$.
\end{Construction}

Next we recall a concrete construction of $\SO(2)$-invariants arising from complexes endowed $\SO(2)$-actions.
Let $\MM$ be a $\kep$-module.
Let $M$ be the object of $\Fun(BS^1,\Mod_k)$ that corresponds
to $\MM$ through $\Fun(BS^1,\Mod_k)\simeq \Comp(\kep)[W^{-1}]$.
We will give an explicit model of the $\SO(2)$-invariants (homotopy fixed points) $M^{\SO(2)}$.
Let $b_{\MM}:\MM[1]\to \MM$ be the map defined by the composite map $\MM[1]\simeq \MM\cdot \epsilon \hookrightarrow \MM\oplus \MM\cdot \epsilon =\MM\otimes_k\kep \to \MM$
where the final map is determined by the module action of $\kep$ on $\MM$.
Since $\epsilon^2=0$, it follows that $b_{\MM}^2=0$. Moreover,
$b_{\MM}d_{\MM}+d_{\MM}b_{\MM}=0$ where $d_{\MM}$ is the differential of $\MM$
(take into account the module action of $\kep$ on $\MM$).
We consider the graded module $\MM[[t]]$ to be $\prod_{i\ge 0}\MM\cdot t^i$.
Here $t$ is a formal variable of cohomological degree $2$.
If $\MM_l$ denotes the part of (homological) degree $l$
of $\MM$,
then the part of (homological) degree $r$ of $\MM[[t]]$
is $\prod_{i\ge 0, l-2i=r}\MM_l\cdot t^i$.
Using $d_{\MM}^2=0$, $b_{\MM}^2=0$,
and $b_{\MM}d_{\MM}+d_{\MM}b_{\MM}=0$,
we see that $d_{\MM}+tb_{\MM}:\MM[[t]]\to \MM[[t]]$ is a square zero map.
We will think of $\MM[[t]]$ as the complex
endowed with the differential $d_{\MM}+tb_{\MM}$.
We observe that $\MM[[t]]$ is an explicit model of
the $\SO(2)$-invariants $M^{\SO(2)}$.
To this end, we choose a cofibrant replacement $K$ of $k$:
\[
\cdots \stackrel{\cdot \epsilon}{\to} k\epsilon \stackrel{0}{\to} k \stackrel{\cdot \epsilon}{\to} k\epsilon \stackrel{0}{\to} k
\]
with respect to the projective model structure.
Then the graded module of Hom complex $\Hom_{\kep}(K,\MM)$
can naturally be identified with $\MM[[t]]$ (if
we denote by $1_{2i}$ the unit $1\in k$ placed in (homological) degree
$2i$, $f:K\to \MM$ corresponds
to $\Sigma_{i=0}^{\infty}f(1_{2i})t^{i}$ in $\MM[[t]]$).
Unwinding the definition we see that 
the differential of the Hom complex corresponds to
$d_{\MM}+tb_{\MM}$.
It follows that $\MM[[t]]$ represents the Hom complex
$\Hom_{\Fun(BS^1,\Mod_k)}(k,M)\simeq M^{\SO(2)}$
where $k$ in $\Fun(BS^1,\Mod_k)$ is the complex
$k$ endowed with the trivial $\SO(2)$-action.
Let $k[t]$ be the free commutative dg algebra generated by $t$
where the (cohomological) degree of $t$ is $2$. The differential
of $k[t]$ is zero.
The Hom complex $\Hom_{\kep}(K,k)$
can be identified with $k[t]$. This means that 
$k[t]=H^*(BS^1,k)$ is a model of $k^{\SO(2)}$.
Let $s:k\to K$ be a canonical section of the resolution $K\to k$.
Then the obvious $k[t]$-module structure on $\MM[[t]]$
corresponds to
\[
\Hom_{\kep}(K,k)\otimes \Hom_{\kep}(K,\MM)\to \Hom_{\kep}(K,\MM),\ \ \phi\otimes f\mapsto f \circ s\circ \phi.
\]
The assignment $\MM\mapsto \MM[[t]]$ defines a functor
\[
\mathbf{F}:\Comp(\kep) \longrightarrow \Comp(k[t])
\]
which preserves quasi-isomorphisms (because every object in $\Comp(\kep)$
is fibrant with respect to the projective model structure).

\begin{Construction}
\label{Finduce}
We apply $\mathbf{F}$ to the diagram~(\ref{diagramB}).
Then we have
\[
\xymatrix{
\GGG_{\CCC}[[t]] \ar[r]^(0.3){\mathbb{A}_{\CCC}^L[[t]]} \ar[d] & \End(HH_\ast(\CCC/k))[[t]] \\
(\GGG_{\CCC}\otimes_kk[\eta])[[t]] \ar[ur]_{\widehat{\mathbb{A}}_{\CCC}^L[[t]]} & 
}
\]
in $\Comp(k[t])$.
We set $(\GGG_{\CCC}\otimes_kk[\eta])((t)):=(\GGG_{\CCC}\otimes_kk[\eta])[[t]][\frac{1}{t}]$ and $\GGG_{\CCC}((t)):=\GGG_{\CCC}[[t]][\frac{1}{t}]$.
Similarly, we write $\End(HH_\ast(\CCC/k))((t))$ for $\End(HH_\ast(\CCC/k))[[t]][\frac{1}{t}]$.
By inverting $t$,
this diagram is extended to the commutative diagram
\[
\label{diagramC}
\vcenter{
\xymatrix{
\GGG_{\CCC}[[t]] \ar[r] \ar[d] & \End(HH_\ast(\CCC/k))[[t]] \ar[r]^i & \End(HH_\ast(\CCC/k))((t)) \\
(\GGG_{\CCC}\otimes_kk[\eta])[[t]] \ar[ur] \ar[r] & (\GGG_{\CCC}\otimes_kk[\eta])((t)). \ar[ur]_{\widehat{\mathbb{A}}_{\CCC}^L((t))} & 
}\tag{3.3}}
\]
\end{Construction}

\begin{Lemma}
\label{mysteriousisom}
We write $k[t^{\pm1}]$ for $k[t][\frac{1}{t}]$.
There exists an isomorphism of dg $k[t^{\pm1}]$-modules
\[
(\GGG_{\CCC}\otimes_kk[\eta])((t))\simeq \GGG_{\CCC}((t))\oplus \GGG_{\CCC}((t))[1]
\]
where the right-hand side is the mapping cone of the identity map of 
$\GGG_{\CCC}((t))$.
\end{Lemma}

\Proof
We first observe that 
$b_{\GGG_{\CCC}}:\GGG_{\CCC}\otimes_kk[\eta][1]\to \GGG_{\CCC}\otimes_kk[\eta]$ is the $k$-linear map determined by
$b_{\GGG_{\CCC}}(g\otimes 1)=0$ and $b_{\GGG_{\CCC}}(g\otimes\eta)=g\otimes 1$.
Thus, 
$d_{\GGG_{\CCC}}+tb_{\GGG_{\CCC}}:(\GGG_{\CCC}\otimes_kk[\eta])[[t]]\to (\GGG_{\CCC}\otimes_kk[\eta])[[t]]$
is given by $(d_{\GGG_{\CCC}}+tb_{\GGG_{\CCC}})(g\otimes 1\cdot t^n)=d_{\GGG_{\CCC}}(g)\otimes 1\cdot t^n$
and $(d_{\GGG_{\CCC}}+tb_{\GGG_{\CCC}})(g\otimes \eta\cdot t^n)=-d_{\GGG_{\CCC}}(g)\otimes \eta\cdot t^n+g\otimes 1\cdot t^{n+1}$.
Here $g\otimes \eta\cdot t^n$ and $g\otimes 1\cdot t^n$
are elements of $(\GGG_{\CCC}\otimes_kk[\eta])[[t]]$
such that
$g\in \GGG_{\CCC}$.
There exists an isomorphism of graded vector spaces
\[
\GGG_{\CCC}((t))\oplus (\GGG_{\CCC}\otimes \eta)((t))\simeq \GGG_{\CCC}((t))\oplus \GGG_{\CCC}((t))[1]
\]
which sends $(g\otimes 1\cdot t^n, g\otimes \eta\cdot t^m)$
to $(g\cdot t^n, g\cdot t^{m+1}[1])$, where $g\cdot t^{m+1}[1]$
indicates an element of $\GGG_{\CCC}((t))[1]$.
The isomorphism transfer the differential $d_{\GGG_{\CCC}}+tb_{\GGG_{\CCC}}$
to the differential $\partial$ of $\GGG_{\CCC}((t))\oplus \GGG_{\CCC}((t))[1]$.
By the computation,
$(\GGG_{\CCC}((t))\oplus \GGG_{\CCC}((t))[1],\partial)$
is the (standard) mapping cone of the identity map
$\GGG_{\CCC}((t))\to \GGG_{\CCC}((t))$.
\QED

The following is standard:

\begin{Proposition}
\label{standard}
Let $f:P\to Q$ be a morphism in $\Comp(k)$, that is, a chain map. 
Let $K=P\oplus Q[-1]$ be the standard mapping cocone of $f$,
that has the natural morphism $\textup{pr}_1:K\to P$ defined as the first projection.
Let $g:C\to P$ be a chain map.
Let $\sigma :C\to Q[-1]$ be a chain homotopy
such that
$d_{Q}\sigma+\sigma d_C=f\circ g$.
Then $g\times \sigma:C\to K$ is a chain map extending $g$.
Moreover, if we consider $\Comp(k)\to \Mod_k=\Comp(k)[W^{-1}]$
and the kernel $\Ker(f)$ of $f$ in $\Mod_k$, then
$g\times \sigma:C\to K$ is a model of the morphism $C\to \Ker(f)$
determined  by $g$ and the homotopy from $f\circ g$ to $0$ associated to $\sigma$
in $\Mod_k$.
\end{Proposition}

\Proof
We prove that
$g\times \sigma$ is compatible with differentials on $C$
and $K$.
Recall that the differential of the mapping cocone
$K=P\oplus Q[-1]$ is defined
by
$P_n\oplus Q_{n+1}\to P_{n-1}\oplus Q_{n}$
which sends $(a,b)$ to $(d_{P}(a), -d_{Q}(b)+f(a))$,
where $P_l$ and $Q_l$ indicate the parts of degree $l$, and
$d_{P}$ and $d_{Q}$ are internal differentials
of $P$ and $Q$,
respectively.
We need to prove
\[
(g(d_{C}(c)), \sigma(d_{C}(c)))=(d_{P}(g(c)), -d_{Q}(\sigma(c))+f(g(c)))
\]
for $c\in C$.
This equality follows from $d_{Q}\sigma+\sigma d_C=f\circ g$ and the fact that
$g$ is a map of complexes.

We prove the latter assertion.
Note that the image of $f \circ \textup{pr}_1:K\to P\to Q$ under $\Comp(k)\to \Mod_k$
and the chain homotopy 
defined as the second projection $K\to Q[-1]$ defines
the pullback square in $\Mod_k$:
\[
\xymatrix{
\Ker(f) \ar[r] \ar[d] & P \ar[d] \\
0 \ar[r] & Q.
}
\]
The space $\Map_{\Mod_k}(C,\Ker(f))\times_{\Map_{\Mod_k}(C,P)}\{g\}$
can be identified with $\Map_{\Map_{\Mod_k}(C,Q)}(f\circ g,0)$ in the natural way.
Thus, 
$g\times \sigma:C\to K$ is a model of the morphism $C\to \Ker(f)$
determined  by the homotopy from $f\circ g$ to $0$ defined by the chain homotopy
$\sigma$.
\QED

\begin{Remark}
By the above proof, an extention $g':C\to K$ such that $\textup{pr}_1\circ g'=g$
bijectively corresponds to a chain homotopy $\sigma$.
\end{Remark}

\begin{Lemma}
\label{extensiontokernel}
Let $i:\End(HH_*(\CCC/k))[[t]]  \to \End(HH_*(\CCC/k))((t))$
be the obvious inclusion and
let $\End(HH_*(\CCC/k))[[t]]\oplus \End(HH_*(\CCC/k))((t))[-1]$
be the mapping cocone of $i$.
Then there exists a morphism $\mathbb{X}$ of complexes
of $k[t]$-modules which makes the diagram
\[
\xymatrix{
  &  \GGG_{\CCC}[[t]] \ar[d]^{\mathbb{A}_{\CCC}^L[[t]]} \ar[ld]_{\mathbb{X}} \\
\End(HH_*(\CCC/k))[[t]]\oplus \End(HH_*(\CCC/k))((t))[-1]\ar[r]^(0.7){\textup{pr}_1}  & \End(HH_*(\CCC/k))[[t]]  \ar[d]^i \\
  &  \End(HH_*(\CCC/k))((t))
}
\]
commute.
\end{Lemma}

\Proof
Based on Proposition~\ref{standard},
we define a map $\mathbb{X}$ of complexes by giving a chain homotopy
$\GGG_{\CCC}[[t]]\to \End(HH_\ast(\CCC/k))((t))[-1]$.
For this purpose, we use a canonical chain homotopy 
$h_{\textup{can}}:\GGG_{\CCC}((t))\oplus \GGG_{\CCC}((t))[1]\to (\GGG_{\CCC}((t))\oplus \GGG_{\CCC}((t))[1])[-1]$ given by $(a,b)\mapsto (0,a)$.
This defines a homotopy from the identity map to the zero map.
Let $J:\GGG_{\CCC}[[t]]\to \End(HH_\ast(\CCC/k))((t))[-1]$ be the composite
\begin{eqnarray*}
\GGG_{\CCC}[[t]]\oplus 0 \hookrightarrow \GGG_{\CCC}((t))\oplus \GGG_{\CCC}((t))[1]\stackrel{h_{\textup{can}}}{\longrightarrow} (\GGG_{\CCC}((t))\oplus \GGG_{\CCC}((t))[1])[-1] &\simeq& (\GGG_{\CCC}\otimes_kk[\eta])((t))[-1]  \\ 
&\stackrel{\widehat{\mathbb{A}}_{\CCC}^L((t))[-1]}{\longrightarrow}& \End(HH_\ast(\CCC/k))((t))[-1].
\end{eqnarray*}
The isomorphism on the right-hand side comes from Lemma~\ref{mysteriousisom}.
We define
$\mathbb{X}$
to be
$\mathbb{A}_{\CCC}^L[[t]]\times J$.
By Proposition~\ref{standard},
$\mathbb{X}$ is a required morphism.
\QED

Consider
the quasi-isomorphism
\[
\End(HH_*(\CCC/k))[[t]]\oplus \End(HH_*(\CCC/k))((t))[-1] \to (\End(HH_*(\CCC/k))((t))/\End(HH_*(\CCC/k))[[t]])[-1]
\]
given by $(a,b)\mapsto\ b\ \textup{mod}\ \End(HH_*(\CCC/k))[[t]]$ (the differential on the right-hand side is zero).
 By Lemma~\ref{extensiontokernel}, we have

\begin{Corollary}
\label{explicitformula}
The composite
\begin{eqnarray*}
\GGG_{\CCC}[[t]] &\stackrel{\mathbb{X}}{\to}& \End(HH_*(\CCC/k))[[t]]\oplus \End(HH_*(\CCC/k))((t))[-1] \\  &\to& (\End(HH_*(\CCC/k))((t))/\End(HH_*(\CCC/k))[[t]])[-1]
\end{eqnarray*}
is given by $\Sigma_{n=0}^\infty g_nt^n\mapsto J(\Sigma_{n=0}^\infty g_nt^n)\ \textup{mod}\ \End(HH_*(\CCC/k))[[t]]$.

\end{Corollary}


\subsection{}

Let
\[
\iota:\End^L(\HH_\bullet(\CCC/k))[[t]]=\End^L(\HH_\bullet(\CCC/k))^{\SO(2)}\to \End^L(\HH_\bullet(\CCC/k))((t))=\End^L(\HH_\bullet(\CCC/k))[[t]][\frac{1}{t}]
\]
be the canonical morphism given by inverting $t$.
We think of $i$ as a morphism of dg Lie algebras. 
The morphism $\iota$ can be regarded as the image of $i$ under $\textup{Lie}^{dg}(k)\to \Lie_k$.
Let $\Ker(\iota)$ be the fiber of $\iota$
in $\Lie_k$.
Since the forgetful functor $\Lie_k \to \Mod_k$ preserves limits, thus
the underlying compelx of $\Ker(\iota)$ is quasi-isomorphic to
the mapping cocone of $i$, that is,
\[
\End(HH_*(\CCC/k))[[t]]\oplus \End(HH_*(\CCC/k))((t))[-1] \to (\End(HH_*(\CCC/k))((t))/\End(HH_*(\CCC/k))[[t]])[-1]
\]
(see Lemma~\ref{extensiontokernel} for the notation).
We will define $\GG_{\CCC}\to \Ker(\iota)$ in $\Lie_k$.

\vspace{2mm}

$\bullet$ In the rest of this section we further assume that $\HH^\bullet(\CCC/k)$
and $\HH_\bullet(\CCC/k)$ belong to $\Perf_k$. Namely,
they are bounded with finite dimensional cohomology.

\vspace{1mm}

We consider the following sequence of symmetric monoidal functors
\[
\Fun(BS^1,\Perf_k)\stackrel{\sim}{\longrightarrow} \Perf_{k[t]}\to \Perf_{k[t^{\pm1}]}
\]
(see Section~\ref{Pr3} for $\Perf$).
The symmetric monoidal functor on the left-hand side
is given by taking $\SO(2)$-invariants. This functor is an equivalence
(see e.g. \cite[Proposition 3.1.4]{Pre}).
The symmetric monoidal functor on the right-hand side
is given by the base change along $k[t]\to k[t^{\pm1}]$.
It gives rise to $\Fun(BS^1,\Lie(\Perf_k))\simeq \Lie(\Fun(BS^1,\Perf_k))\stackrel{\sim}{\to} \Lie(\Perf_{k[t]}) \to \Lie(\Perf_{k[t^{\pm1}]})$.
Let 
\[
\xymatrix{
\GG_{\CCC}((t)) \ar[r] \ar[dr]_{A^L_{\CCC}((t))} & (\GG_{\CCC}^{S^1})((t)) \ar[d] \\
  &  \End(\HH_\bullet(\CCC/k))((t)) 
}
\]
be the diagram in $\Lie(\Perf_{k[t^{\pm1}]})$ which is defined as
the image of the diagram~(\ref{diagramA}).

\begin{Lemma}
\label{null}
$(\GG_{\CCC}^{S^1})((t))$ is equivalent to $0$ in $\Lie(\Perf_{k[t^{\pm1}]})$.
\end{Lemma}

\Proof
It will suffice to prove that
the underlying complex $(\GG_{\CCC}^{S^1})((t))$ is quasi-isomorphic to $0$.
Note that
$(\GG_{\CCC}^{S^1})^{\SO(2)}$ is equivalent to $\GG_{\CCC}\simeq \GG_{\CCC}^{S^1/\SO(2)}$
endowed with the $k[t]$-module structure induced by the restriction along
$k[t]=k^{\SO(2)}\to (k^{S^1})^{\SO(2)}=k$.
Thus, $(\GG_{\CCC}^{S^1})((t))=(\GG_{\CCC}^{S^1})^{\SO(2)}\otimes_{k[t]}k[t^{\pm1}]\simeq 0$.
\QED

\begin{Construction}
\label{liftconstruction}
We will define $\chi^L:\GG_{\CCC}\to \Ker(\iota)$.
By Lemma~\ref{null}
there exists an (essentially unique)
homotopy $h$ from $\textup{id}_{(\GG_{\CCC}^{S^1})((t))}$
to $0$, that is an equivalence in
$\Map_{\Lie(\Perf_{k[t^{\pm1}]})}((\GG_{\CCC}^{S^1})((t)),(\GG_{\CCC}^{S^1})((t)))$.
Composing $h$ with $\GG_{\CCC}((t)) \to (\GG_{\CCC}^{S^1})((t))$
and $(\GG_{\CCC}^{S^1})((t)) \to \End(\HH_\bullet(\CCC/k))((t))$
we obtain
a homotopy $h_0$ from $A^L_{\CCC}((t))$ to $0$.
Consider the forgetful functor $\Lie(\Perf_{k[t]})\to \Lie_k$.
Since 
$\GG_{\CCC}[[t]]\to \End(\HH_\bullet(\CCC/k))((t))$
factors through $A^L_{\CCC}((t)):\GG_{\CCC}((t))\to \End(\HH_\bullet(\CCC/k))((t))$,
the image of $h_0$ induces a homotopy $h_1$ from the composite 
\[
\GG_{\CCC}\to \GG_{\CCC}[[t]]\to \End(\HH_\bullet(\CCC/k))[[t]]\to \End(\HH_\bullet(\CCC/k))((t))
\]
to $0$ in $\Lie_k$.
It gives rise to $\chi^L:\GG_{\CCC}\to \Ker(\iota)$ which makes the diagram
\[
\xymatrix{
\GG_{\CCC} \ar[rd]^{\chi^L} \ar[rdd] \ar[rrd] &     &    \\
                   & \Ker(\iota) \ar[r] \ar[d] &  \End(\HH_\bullet(\CCC/k))[[t]] \ar[d]^\iota \\
                   & 0               \ar[r] &  \End(\HH_\bullet(\CCC/k))((t))
}
\]
commute in $\Lie_k$.
\end{Construction}

\begin{Lemma}
\label{underlyingmodel}
The underlying map $\chi:\GG_{\CCC} \to \Ker(\iota)$ of $\chi^L$ in $\Mod_k$
is represented by $\mathbb{X}_0:\GGG_{\CCC}\to \GGG_{\CCC}[[t]]\stackrel{\mathbb{X}}{\to} \End(HH_*(\CCC/k))[[t]]\oplus \End(HH_*(\CCC/k))((t))[-1]$.
\end{Lemma}

\Proof
The underlying map $\chi:\GG_{\CCC} \to \Ker(\iota)$ of $\chi^L$ in $\Mod_k$
is induced by the underlying homotopy $\bar{h}_1$ of $h_1$ in $\Mod_k$
(see Construction~\ref{liftconstruction} for $h_1$).
The diagram~(\ref{diagramB}) is a model of the diagram~(\ref{diagramA}) at the level of $\Fun(BS^1,\Mod_k)$,
and the diagram~(\ref{diagramC}) in Construction~\ref{Finduce}
is a model of
\[
\xymatrix{
\GG_{\CCC}[[t]] \ar[r] \ar[d] & \End^L(\HH_\bullet(\CCC/k))[[t]] \ar[r] & \End^L(\HH_\bullet(\CCC/k))((t)) \\
(\GG_{\CCC}^{S^1})[[t]] \ar[ur] \ar[r] & (\GG_{\CCC}^{S^1})((t)). \ar[ur] & 
}
\]
By Lemma~\ref{mysteriousisom},
$(\GGG_{\CCC}\otimes k[\eta])((t))$ is isomorphic to the standard mapping cocone
$\GGG_{\CCC}((t))\oplus \GGG_{\CCC}((t))[1]$.
There exists an essentially unique homotopy from the identity map 
of $\GGG_{\CCC}((t))\oplus \GGG_{\CCC}((t))[1]$
to the zero map in $\Mod_k$. In Lemma~\ref{extensiontokernel},
we use the canonical chain homotopy
$h_{\textup{can}}$ to define the chain homotopy $J$ that determines $\mathbb{X}$
(cf. Proposition~\ref{standard}).
Thus, by comparing the definition of $\bar{h}_1$
and the definition of $\mathbb{X}$, we see that 
$\chi$
is represented by $\mathbb{X}_0:\GGG_{\CCC}\to \GGG_{\CCC}[[t]]\stackrel{\mathbb{X}}{\to} \End(HH_*(\CCC/k))[[t]]\oplus \End(HH_*(\CCC/k))((t))[-1]$.
\QED

\begin{Lemma}
\label{IMLemma}
Consider $(\End(HH_*(\CCC/k))((t))/\End(HH_*(\CCC/k))[[t]])[-1],\ \textup{zero bracket})$ to be a dg Lie algebra endowed with the zero differential and the zero bracket.
Then there exists an equivalence in $\Lie_k$:
\[
\Ker(\iota) \simeq (\End(HH_*(\CCC/k))((t))/\End(HH_*(\CCC/k))[[t]])[-1],\ \textup{zero bracket}).
\]
\end{Lemma}

\Proof
We identify $\iota$ with the image of $i:\End(HH_*(\CCC/k))[[t]]  \to \End(HH_*(\CCC/k))((t))$ under $\Lie^{dg}(k)\to \Lie_k$.
Since $i$ satisfies the assumption of \cite[Proposition 3.4]{IaMa},
it follows from {\it loc.cit.} that
$\Ker(\iota)$ is equivalent to a dg Lie algbera with the zero bracket. Thus the assertion holds.
\QED

\begin{Lemma}
\label{injectivewin}
Suppose that
\[
H^*(\chi^L):H^*(\GG_{\CCC})\to H^*(\Ker(\iota))\simeq (\End(HH_*(\CCC/k))((t))/\End(HH_*(\CCC/k))[[t]])[-1]
\]
is an injective map of graded vector spaces.
Then $\GG_{\CCC}$ is quasi-abelian. Namely, it is equivalent to
a dg Lie algebra with the zero bracket. 
\end{Lemma}

\Proof
Let $V$ be the image of $H^*(\chi^L)$. It is a graded subvector space
which is contained in $E:=(\End(HH_*(\CCC/k))((t))/\End(HH_*(\CCC/k))[[t]])[-1]$.
According to Lemma~\ref{IMLemma}, we can take $E$
endowed with the zero differential and the zero bracket
to be a dg Lie algbera which is a model of $\Ker(\iota)\in \Lie_k$.
We regard $V$ as a dg Lie algebra with the zero differential and
the zero bracket.
There exists a genuine contraction map
\[
c:E=(\End(HH_*(\CCC/k))((t))/\End(HH_*(\CCC/k))[[t]])[-1]\to V
\]
of dg Lie algebras such that $V\to E\to V$ is the identity.
The composition of $c$ and $\chi^L$ defines an equivalence 
$\GG_{\CCC}\simeq (V,\ \textup{zero bracket})$ in $\Lie_k$.
\QED

\begin{Proposition}
\label{mainreduction3}
Suppose that
\[
H^*(\mathbb{X}_0):H^*(\GGG_{\CCC})=H^*(\GG_{\CCC})\to H^*(\Ker(\iota))=
(\End(HH_*(\CCC/k))((t))/\End(HH_*(\CCC/k))[[t]])[-1]
\]
is injective.
Then $\GG_{\CCC}$ is quasi-abelian.
\end{Proposition}

\Proof
Proposition follows from Lemma~\ref{underlyingmodel} and Lemma~\ref{injectivewin}.
\QED

\section{Around Contraction morphisms}

\label{contractionsection}

In this section we will relate two morphisms which we will denote by $I_{\CCC}$
and $\gamma_{\CCC}$.
The morphism $I_{\CCC}$ is called the contraction morphism.
This is an associative algebra map that
has a geometric description in terms of factorization homology
and the $\HH_\bullet(\HH^\bullet(\CCC/k))$-module $\HH_\bullet(\CCC/k)$
(see Definition~\ref{defcontraction} and Remark~\ref{remcontraction}).
The morphism 
$\gamma_{\CCC}$ is a dg Lie algebra map which is directly related to the morphism in Propostion~\ref{mainreduction3}
(see Lemma~\ref{representkernel}).
We prove Proposition~\ref{desiredcompatibility}
which identifies $I_{\CCC}$ with $\gamma_{\CCC}$ as a morphism in $\Mod_k$.
As a result we obtain Proposition~\ref{ibunnot}.

\subsection{Factorization Homology}
We first review the factorization homology.
The review is minimal and quick because the factorization homology will be applied
only to the Euclidean surface $\RRR^2$
and the cylinder $\RRR\times S^1$. 
We briefly recall the definition of factorization homology. 
We refer the reader to \cite{AF}, \cite{Ho} for the theory of factorization homology
we need (in \cite{HA},
the theory is developed under the name of topological chiral homology).
Let $(\Mfldf_n)^{\otimes}$ be the symmetric monoidal $\infty$-category of framed smooth $n$-manifolds
such that the mapping spaces are the space
of smooth embeddings endowed the data of compatibility of framings
(see e.g. \cite[2.1]{AF} or \cite{Ho} for the detailed account).
The symmetric monoidal structure is given by disjoint union.
We write $\Mfldf_n$ for the underlying $\infty$-category.
In a nutshell, the factorization homology is defined as the functor
\[
\int_{(-)}(\bullet):\Alg_n(\Mod_k)\times \Mfldf_n \longrightarrow \Mod_k
\]
which carries $(B,M)$ to $\int_MB$. 
We refer to $\int_MB$ as the factorization homology
of $M$ in coefficients in $B$.
The strategy of the construction of $\int_{(-)}(\bullet)$ is as follows.
Let $\Diskf_n$ be the full subcategory of  $\Mfldf_n$ spanned by manifolds
which are diffeomorphic to a (possibly empty) finite disjoint union of $\RRR^n$.
Let $p:(\Diskf_n)^{\otimes}\to (\Mfldf_n)^{\otimes}$ denote the fully faithful symmetric monoidal functor.
Let $\Fun^\otimes((\Diskf_n)^{\otimes},\Mod_k^\otimes)$ 
and $\Fun^\otimes((\Mfldf_n)^{\otimes},\Mod_k^\otimes)$
denote the $\infty$-categories
of symmetric monoidal functors.
There exists a canonical equivalence $\Fun^\otimes((\Diskf_n)^{\otimes},\Mod_k^\otimes)\simeq \Alg_n(\Mod_k)$
(obtained by identifying $(\Diskf_n)^{\otimes}$ with a symmetric monoidal envelope of the operad $\eenu$,
see \cite[Section 2.2.4]{HA} for monoidal evelopes).
Consider the following adjoint pair
\[
\xymatrix{
p_!: \Fun^\otimes((\Diskf_n)^{\otimes},\Mod_k^\otimes)  \ar@<0.5ex>[r] &  \Fun^\otimes((\Mfldf_n)^{\otimes},\Mod_k^\otimes)  : p^*   \ar@<0.5ex>[l]  
}
\]
where $p^*$ is determined by composition with $p$.
The left adjoint $p_!$ sends $\beta:(\Diskf_n)^{\otimes}\to \Mod_k^\otimes$
to a symmetric monoidal (operadic) left Kan extension $(\Mfldf_n)^{\otimes}\to \Mod_k^\otimes$ of $\beta$ along $p$.
We consider the composite functor  
$\Alg_n(\Mod_k)\simeq \Fun^\otimes((\Diskf_n)^{\otimes},\Mod_k^\otimes)\to  \Fun^\otimes((\Mfldf_n)^{\otimes},\Mod_k^\otimes)\to \Fun(\Mfldf_n,\Mod_k)$
where the final functor sends symmetric monoidal functors
to underlying functors (obtained by forgetting symmetric monoidal structure).
Passing to the adjoint of the composite functor, we have $\int_{(-)}(\bullet)$.
By the definition of $p_!$, there exists a canonical equivalence
$\int_{\RRR^n}B\simeq B$.
If $M=\RRR\times S^1$ and $B\in \Alg_2(\Mod_k)$, $\int_MB\simeq \HH_\bullet(B/k)$ (see e.g. \cite[Proposition 7.11]{ID}).
Let us regard $\beta:(\Diskf_n)^{\otimes}\to \Mod_k^\otimes$ as an $\eenu$-algebra $B$ in $\Mod_k$.
For any $M\in \Mfldf_n$, 
$\int_{M}B$ is the image of $M$ under $p_!(\beta)$.
The factorization homology $\int_{M}B$ can naturally be identified with a colimit of 
$(\Diskf_n)_{/M}\to \Diskf_n\to \Mod_k$
where the second functor
is the {\it underlying functor}
of the symmetric monoidal functor $\beta:(\Diskf_n)^\otimes\to \Mod_k^\otimes$.

\subsection{Contraction morphisms}
Consider an open embedding $j:\RRR\to S^1$
which identify $\RRR$ with an open set of $S^1$.
Then it gives rise to
an open embedding
$\textup{id}\times j:\RRR\times \RRR\to \RRR \times S^1$ and
the induced morphism
\[
\HH^\bullet (\CCC/k) \simeq \int_{\RRR^2}\HH^\bullet(\CCC/k)\to \int_{\RRR \times S^1}\HH^\bullet(\CCC/k)\simeq \HH_\bullet(\HH^\bullet(\CCC/k)).
\]

Let 
$A_{\CCC}:\HH_\bullet(\HH^\bullet(\CCC/k)/k)\to \End(\HH_\bullet(\CCC/k))$
denote the morphism of associative algebras (morphism in $\Alg_1(\Mod_k)$
determined by the action of $\HH_\bullet(\HH^\bullet(\CCC/k)/k)$
on $\HH_\bullet(\CCC/k)$ (see Section~\ref{KSalgebrasection}).

\begin{Definition}
\label{defcontraction}
Let $I_{\CCC}$ be the composition
\[
\HH^\bullet (\CCC/k)\to \HH_\bullet(\HH^\bullet(\CCC/k)/k)\stackrel{A_{\CCC}}{\longrightarrow} \End(\HH_\bullet(\CCC/k))
\]
in $\Mod_k$
which we refer to as the contraction morphism.
\end{Definition}

\begin{Remark}
\label{remcontraction}
Consider the open embedding $j_r:\RRR^{\sqcup r}=\RRR\sqcup\ldots \sqcup \RRR\to \RRR$
and the map $j_r\times \textup{id}_{S^1} :\RRR^{\sqcup r}\times S^1\to \RRR\times S^1$.
The covariant fucntoriality of factorization homology
defines $\bigl(\int_{\RRR\times S^1}\HH^\bullet(\CCC/k)\bigr)^{\otimes r+1}\to \int_{\RRR \times S^1}\HH^\bullet(\CCC/k)$.
Moreover, the functoriality of the first factor $\RRR$ promotes $\int_{\RRR \times S^1}\HH^\bullet(\CCC/k)$ to an object of $\Alg_1(\Mod_k)$,
and $\int_{\RRR^2}\HH^\bullet(\CCC/k)\to  \int_{\RRR \times S^1}\HH^\bullet(\CCC/k)$
is lifted to a morphism in $\Alg_1(\Mod_k)$.
Consequently, $I_{\CCC}$ is promoted to a morphism in $\Alg_1(\Mod_k)$.
We will not use this fact in this note.
\end{Remark}

\subsection{Comparison with the extended Lie algebra action}
Let $\ast\to S^1$ be a point of $S^1$, and
$\GG_{\CCC}^{S^1}\to  \GG_{\CCC}$ the induced morphism
of dg Lie algebras.
The kernel of $\GG_{\CCC}^{S^1}\to \GG_{\CCC}$
is $\Omega_L(\GG_{\CCC})$.
Consider the composition
\[
\Omega_L(\GG_{\CCC}) \to \GG_{\CCC}^{S^1}\stackrel{\widehat{A}^L_{\CCC}}{\longrightarrow} \End^L(\HH_\bullet(\CCC/k)).
\]
Let $\gamma_{\CCC}$ be the induced morphism
$\HH^\bullet(\CCC/k)\simeq \Omega_L(\GG_{\CCC}) \to \End(\HH_\bullet(\CCC/k))$
in $\Mod_k$ where $\HH^\bullet(\CCC/k)\simeq \Omega_L(\GG_{\CCC})$
is the canonical equivalence in $\Mod_k$ (see the first paragraph of
the proof of Lemma~\ref{identifyI}).

\begin{Proposition}
\label{desiredcompatibility}
The morphism $\gamma_{\CCC}$ is equivalent to $I_{\CCC}$ as a morphism in $\Mod_k$
through the canonical equivalences $\Omega_L(\GG_{\CCC})\simeq \HH^\bullet(\CCC/k)$
and $\End^L(\HH_\bullet(\CCC/k))\simeq \End(\HH_\bullet(\CCC/k))$ in $\Mod_k$.
\end{Proposition}

\begin{Lemma}
\label{middlesquare}
Let $U_1(\Omega_L(\GG_{\CCC}))\to U_1(\GG_{\CCC}^{S^1})$
be the morphism induced by the canonical morphism
$\Omega_L(\GG_{\CCC})\to \GG_{\CCC}^{S^1}$.
Consider $U_2(\GG_{\CCC})\simeq \int_{\RRR^2} U_2(\GG_{\CCC})\to \int_{\RRR\times S^1}U_2(\GG_{\CCC})$
induced by $\textup{id}\times j:\RRR \times \RRR \to \RRR\times S^1$.
Then there exists a diagram in $\Mod_k$
\[
\xymatrix{
U_1(\Omega_L(\GG_{\CCC})) \ar[r]^{\simeq} \ar[d] & U_2(\GG_{\CCC}) \ar[d] \\
U_1(\GG_{\CCC}^{S^1}) \ar[r]^{\simeq} &  \int_{\RRR\times S^1}U_2(\GG_{\CCC})
}
\]
which commutes up to a homotopy such that
the horizontal arrows are equivalences.
\end{Lemma}

Before the proof, we recall the nonabelian Poincar\'e duality which we will use 
in the proof.
Let $M$ be a framed (smooth) $n$-manifold and let $L$ be a dg Lie algebra.
Let $Ch_\bullet:\Lie_k\to \Mod_k$ be the Chevalley-Eilenberg chain functor.
Let $\Map_{c}(M,L)$ denote the compactly supported 
cochains of $M$ with coefficients
$L$.
That is, $\Map_{c}(M,L)$ can be identified with
the kernel of $L^{M^+}\to L^{+}\simeq L$ induced by the restriction to 
the added point $+$ where 
$M^+$ is the one-point compactification of $M$, and
$L^{N}$ denotes the cotensor with a space $N$
in $\Lie_k$.
The nonabelian Pioncar\'e duality for dg Lie algebras \cite[Proposition 5.13]{AF} says that there exists a canonically defined equivalence
\[
\int_MCh_\bullet(\Omega_L^n(L))\simeq Ch_\bullet(\Map_c(M,L))
\] which is functorial
with respect to embeddings of framed $n$-manifolds.
Recall from Section~\ref{Pre2} that $U_1\simeq Ch_\bullet\circ \Omega_L$
and $U_n\simeq U_1\circ \Omega_L^{n-1}$ so that 
$U_n\simeq Ch_\bullet\circ \Omega_L^{n}$.
It follows that $\int_MU_n(L)\simeq  Ch_\bullet(\Map_c(M,L))$.

{\it Proof of Lemma~\ref{middlesquare}.}
Applying the above equivalence to $\RRR \times \RRR \to \RRR\times S^1$
and $\GG_{\CCC}$,
we obtain the diagram
\[
\xymatrix{
Ch_\bullet(\Map_c(\RRR^2,\GG_{\CCC})) \ar[r]^(0.6){\simeq} \ar[d] & \int_{\RRR^2}U_2(\GG_{\CCC}) \ar[d] \\
Ch_\bullet(\Map_c(\RRR\times S^1,\GG_{\CCC})) \ar[r]^(0.6){\simeq}  & \int_{\RRR\times S^1}U_2(\GG_{\CCC}) 
}
\]
which commute up to a canonical homotopy.
Moreover, by using equivalences
$\Map_c(\RRR,L)\simeq \Omega_L(L)$
and $\Map_c(S^1,L)=L^{S^1}$ for $L\in \Lie_k$, we see that 
the left vertical arrow can be identified with
\[
Ch_\bullet(\Map_c(\RRR,\Map_c(\RRR,\GG_{\CCC})))\simeq Ch_\bullet(\Omega_L(\Omega_L(\GG_{\CCC}))) \to Ch_\bullet(\Map_c(\RRR,\Map_c(S^1,\GG_{\CCC})))=Ch_\bullet(\Omega_L(\GG_{\CCC}^{S^1}))
\]
where the morphism is induced by 
$\Map_c(\RRR,\GG_{\CCC})\to \Map_c(S^1,\GG_{\CCC})$
that is identified with $\Omega_L(\GG_{\CCC})=0\times_{\GG_{\CCC}}\GG_{\CCC}^{S^1}\to \GG_{\CCC}^{S^1}$ arising from
 $\GG_{\CCC}^{S^1}\to \GG_{\CCC}$.
By the equivalence $U_n\simeq Ch_\bullet\circ \Omega_L^{n}$ 
between functors $\Lie_k\to \Alg_n(\Mod_k)$,
it can be identified with
$U_2(\GG_{\CCC})\simeq U_1(\Omega_L(\GG_{\CCC}))\to U_1(\GG_{\CCC}^{S^1})$ induced by $\Omega_L(\GG_{\CCC})\to \GG_{\CCC}^{S^1}$.
The assertion follows.
\QED

\begin{Lemma}
\label{identifyI}
There exists a diagram in $\Mod_k$
\[
\xymatrix{
\Omega_L(\GG_{\CCC})\ar[r]^{\simeq} \ar[d]^s & \HH^\bullet(\CCC/k) \ar[d]  \ar[rd]^{\textup{id}}& \\
U_1(\Omega_L(\GG_{\CCC})) \ar[r]^{\simeq} & U_2(\GG_{\CCC}) \ar[r] & \HH^\bullet(\CCC/k)
}
\]
which commutes up to a canonical homotopy.
Here the left vertical arrow $s$ is the canonical unit map arising from the adjoint
pair $U_1:\Lie_k\rightleftarrows \Alg_1(\Mod_k):res_{\eone/Lie}$.
The morphism $U_2(\GG_{\CCC})\to \HH^\bullet(\CCC/k)$
is the counit map of the adjoint pair $U_2:\Lie_k\rightleftarrows  \Alg_2(\Mod_k):res_{\etwo/Lie}$.
\end{Lemma}

{\it Proof of Lemma~\ref{identifyI}.}
We use the notation in Lemma~\ref{identifyI2}.
We first note that $\Omega_L(\GG_{\CCC})$ is defined 
to be $\Omega_L\circ B_L\circ \Alg_1(res_{\eone/Lie})(\HH^\bullet(\CCC/k))$
which is canoncially equivalent to $\Alg_1(res_{\eone/Lie})(\HH^\bullet(\CCC/k))$
(see (2) of Lemma~\ref{identifyI2}). 
Since the composite
$\Alg_2(\Mod_k)\simeq \Alg_1(\Alg_1(\Mod_k))\stackrel{\Alg_1(res_{\eone/Lie})}{\longrightarrow} \Alg_1(\Lie_k^{\times})\stackrel{\textup{forget}}{\longrightarrow} \Mod_k$ is equivalent to the evident forgetful functor,
it follows that there exists
a canonical equivalence 
$\Omega_L(\GG_{\CCC})\simeq \HH^\bullet(\CCC/k)$
in $\Mod_k$.

We will prove that the following diagram commutes up to a homotopy:
\[
\xymatrix{
 \Omega_L(\GG_{\CCC}) \ar[rd] \ar[r]^(0.4){u[-1]} & \Omega_L(\textup{Free}_{Lie}(\GG_{\CCC})) \ar[r]^v & U_1(\Omega_L(\textup{Free}_{Lie}(\GG_{\CCC}))) \ar[r] \ar[d]^\simeq& U_1(\Omega_L(\GG_{\CCC})). \ar[d]^\simeq& \\
 & \Free_{\etwo}(\Omega_L(\GG_{\CCC})) \ar[r]^{\simeq} \ar[ru]^\simeq & U_2(\textup{Free}_{Lie}(\GG_{\CCC})) \ar[r]  & U_2(\GG_{\CCC})  \ar[r] & \HH^\bullet(\CCC/k)
}
\]
The natural equivalence
$U_2\simeq U_1\circ \Omega_L:\Lie_k\to \Alg_2(\Mod_k)$
applied to the canonical morphism
$\textup{Free}_{Lie}(\GG_{\CCC})\to \GG_{\CCC}$
gives the right square diagram which commutes up to a canonical homotopy.
Taking into account Lemma~\ref{identifyI2}, we have the middle triangle which commutes 
up to a canonical homotopy (compare the compositions of left adjoints).
Next we will show that the left triangle commutes.
The unit map of the adjoint pair
$(\Omega_L\circ \textup{Free}_{Lie}, res_{Lie/M}\circ B_L)$
induces
$u:\GG_{\CCC}\to \Omega_L(\textup{Free}_{Lie}(\GG_{\CCC}))[1]$ in $\Mod_k$.
The unit map of the adjoint pair
$(\Alg_1(U_1), \Alg_1(res_{\eone/Lie}))$
induces 
$v:\Omega_L(\textup{Free}_{Lie}(\GG_{\CCC})) \to U_1(\Omega_L(\textup{Free}_{Lie}(\GG_{\CCC})))$.
The unit map of the composition of the adjoint pairs in Lemma~\ref{identifyI2}
induces $v[1]\circ u: \GG_{\CCC}\to \Omega_L(\textup{Free}_{Lie}(\GG_{\CCC}))[1]\to U_1(\Omega_L(\textup{Free}_{Lie}(\GG_{\CCC})))[1]$.
Applying the shift functor $[-1]$ we have $v\circ (u[-1])$.
The morphism $\Omega_L(\GG_{\CCC})\to \textup{Free}_{\etwo}(\Omega_L(\GG_{\CCC}))$ is the morphism induced by the
unit map of the adjoint pair $\textup{Free}_{\etwo}:\Mod_k\rightleftarrows \Alg_2(\Mod_k):res_{\etwo/M}$.
Note that the composite functor
$[-1]\circ res_{Lie/M}\circ B_L\circ \Alg_1(res_{\eone/Lie})$ is equivalent to 
the forgetful functor $res_{\etwo/M}$
so that $\Alg_1(U_1)\circ\Omega_L\circ \textup{Free}_{Lie}\circ  [1]\simeq \textup{Free}_{\etwo}$ where $[1]$ and $[-1]$ indicate the shift functors.
It follows that $\Omega_L(\GG_{\CCC})\to \textup{Free}_{\etwo}(\Omega_L(\GG_{\CCC}))$ determined by $(\textup{Free}_{\etwo},res_{\etwo/M})$
is canonically equivalent to $v\circ (u[-1])$. Namely, the left triangle
commutes up to a homotopy.
Taking into account the sequence of adjoint pairs in Lemma~\ref{identifyI2}
we observe that the composite $\Omega_L(\GG_{\CCC})\to U_1(\Omega_L(\GG_{\CCC}))$ is $s$.
It remains to prove that
the composite $\Omega_L(\GG_{\CCC})\to \textup{Free}_{\etwo}(\Omega_L(\GG_{\CCC}))\to \HH^\bullet(\CCC/k)$ gives the canonical
equivalence $\Omega_L(\GG_{\CCC})\simeq \HH^\bullet(\CCC/k)$.
To this end, consider the sequence of adjoint pairs
\[
\xymatrix{
\Mod_k \ar@<0.5ex>[r]^{[1]}  & \Mod_k \ar@<0.5ex>[r]^{\textup{Free}_{Lie}} \ar@<0.5ex>[l]^{[-1]}  & \Lie_k \ar@<0.5ex>[l]^{res_{Lie/M}}  \ar@<0.5ex>[r]^{U_2} & \Alg_2(\Mod_k) \ar@<0.5ex>[l]^{res_{\etwo/Lie}}
}
\]
such that $[-1]\circ res_{Lie/M}\circ res_{\etwo/Lie}=res_{\etwo/M}$,
and $U_2\circ \textup{Free}_{Lie}\circ [1]\simeq \textup{Free}_{\etwo}$.
For ease of notation, 
we set $F=\textup{Free}_{Lie}$, $F'=U_2$, $R=res_{Lie/M}$, and $R'=res_{\etwo/Lie}$.
Unfolding the definition, the composite is given by the composite of the sequence
\begin{eqnarray*}
\Omega_L(\GG_{\CCC})=\GG_{\CCC}[-1] &=& ([-1]\circ R) (R'(\HH^\bullet(\CCC/k))) \\
&\simeq& ([-1]\circ R\circ R')(\HH^\bullet(\CCC/k)) \\
&\to& ([-1]\circ R\circ R')\circ (F'\circ F\circ [1])\circ ([-1]\circ R\circ R')(\HH^\bullet(\CCC/k)) \\
&\simeq&  [-1]\circ R\circ R'\circ F'\circ (F\circ R)\circ R'(\HH^\bullet(\CCC/k)) \\
&\to&  [-1]\circ R\circ R'\circ (F'\circ R')(\HH^\bullet(\CCC/k)) \\
&\to&  ([-1]\circ R\circ R')(\HH^\bullet(\CCC/k))
\end{eqnarray*}
where the first arrow is induced by the unit map $\textup{id}\to ([-1]\circ R\circ R')\circ (F'\circ F\circ [1])$, the second arrow is induced by
the counit map $F\circ R\to \textup{id}$, and the third arrow is induced by
the counit map $F'\circ R'\to \textup{id}$.
By the counit-unit relation, the composite from the second line to the final line is
equivalent to the identity map.
Thus, our assertion follows.
\QED

{\it Proof of Proposition~\ref{desiredcompatibility}.}
Consider the diagram
\[
\xymatrix{
& & \Omega_L(\GG_{\CCC})\ar[r]^{\simeq} \ar[d] \ar[dl] & \HH^\bullet(\CCC/k) \ar[d]  \ar[rd]^{\textup{id}}& & \\
\GG_{\CCC} \ar[d] \ar[r] & \GG_{\CCC}^{S^1} \ar[dr] & U_1(\Omega_L(\GG_{\CCC})) \ar[r]^{\simeq} \ar[d]  & U_2(\GG_{\CCC}) \ar[r] \ar[d] & \HH^\bullet(\CCC/k) \ar[d]  \ar[dr]^{I_{\CCC}} &  \\
U_1(\GG_{\CCC}) \ar[rr] & & U_{1}(\GG_{\CCC}^{S^1}) \ar[r]^{\simeq} & \int_{\RRR\times S^1} U_2(\GG_{\CCC}) \ar[r] & \int_{\RRR\times S^1}\HH^\bullet(\CCC/k) \ar[r]^{A_{\CCC}} & \End(\HH_\bullet(\CCC/k))
}
\]
in $\Mod_k$.
The right upper trapezoid comes from Lemma~\ref{identifyI}.
The lower middle square comes from Lemma~\ref{middlesquare}.
Both diagrams commute up to homotopy.
The right square is defined by
applying the fatorization homology $\int_{\RRR^2}(-)\to \int_{\RRR\times S^1}(-)$
to $U_{2}(\GG_{\CCC})\to \HH^\bullet(\CCC/k)$ so that it commutes
up to a canonical homotopy.
By the definition of $I_{\CCC}$, the right triangle commutes.
Note that for $L\in \Lie_k$ the canonical morphism $L\to U_1(L)$
is defined as the map induced by the unit map of the adjoint pair
$(U_1,res_{\eone/Lie})$. It follows that $L\to U_1(L)$ is functorial
with respect to $L\in \Lie_k$. 
The left trapezoid and the left triangle are defined by
the canonical maps so that they commute up to homotopy.

Now we will prove our assertion by using this diagram.
Note that by definition, $\gamma_{\CCC}$ is the composite
$\Omega_L(\GG_{\CCC})\to U_1(\GG_{\CCC}^{S^1})\to\int_{\RRR\times S^1} U_2(\GG_{\CCC}) \to \int_{\RRR\times S^1}\HH^\bullet(\CCC/k) \to \End(\HH_\bullet(\CCC/k))$ in the above diagram.
Since the above diagram commutes up to homotopy, it follows that
$\gamma_{\CCC}$ is equivalent to $I_{\CCC}$ (via the canonical equivalence
$\Omega_L(\GG_{\CCC})\simeq \HH^\bullet(\CCC/k)$).
\QED

\begin{Lemma}
\label{representkernel}
The canoncial map
$\Omega_L(\GG_{\CCC})\to \GG_{\CCC}^{S^1}$ (in $\Mod_k$)
is represented by the canonical inclusion
$\GGG_{\CCC}\otimes_k\eta \hookrightarrow \GGG_{\CCC}\otimes_kk[\eta]$.
\end{Lemma}

\Proof
Consider the map $\ast\to S^1$
and the induced morphism
$\Hom_k(k[\epsilon],\GG_{\CCC})\simeq \GG_{\CCC}^{S^1}\to \GG_{\CCC}^{\ast}=\Hom_k(k,\GG_{\CCC})\simeq \GG_{\CCC}$
where $\Hom_k(-,-)$
is the Hom complex defined as an object of $\Mod_k$.
(we may think that the induced map is determined by the composition with $k=H_*(\ast,k)\to H_*(S^1,k)=k[\epsilon]$).
The morphism $\GGG_{\CCC}\otimes_kk[\eta]\to \GGG_{\CCC}$
induced by $k[\eta]\to k$
represents $\GG_{\CCC}^{S^1}\to \GG_{\CCC}^{\ast}=\GG_{\CCC}$.
In the category $\Comp(k)$ (endowed with the projective model structure),
every object is fibrant, and
 $\GGG_{\CCC}\otimes k[\eta]\to \GGG_{\CCC}$ is a fibration.
Thus, $\GGG_{\CCC}\otimes_k\eta$
(with the canonical inclusion $\GGG_{\CCC}\otimes_k\eta \hookrightarrow \GGG_{\CCC}\otimes_kk[\eta]$) represents a homotopy fiber of
$\GGG_{\CCC}\otimes_kk[\eta]\to \GGG_{\CCC}$.
Namely, the kernel (homotopy fiber) $\Omega_L(\GG_{\CCC})$ of $\GG_{\CCC}^{S^1}\to \GG_{\CCC}^{\ast}=\GG_{\CCC}$ is represented by $\GGG_{\CCC}\otimes_k\eta$.
\QED

\begin{Corollary}
\label{kernel2}
The composite 
$\mathbb{I}_{\CCC}:\GGG_{\CCC}\otimes\eta \to \GGG_{\CCC}\otimes_kk[\eta]\stackrel{\widehat{\mathbb{A}}^L_{\CCC}}{\longrightarrow} \End(HH_*(\CCC/k))$
represents a morphism $\HH^\bullet(\CCC/k)=\Omega_L(\GG_{\CCC}) \to \End(\HH_\bullet(\CCC/k))$
which is equivalent to $I_{\CCC}$ in $\Mod_k$.
See also Construction~\ref{modelconstruction} for $\mathbb{I}_{\CCC}$.
\end{Corollary}

\Proof
Combine Proposition~\ref{desiredcompatibility}
and Lemma~\ref{representkernel}.
\QED

\begin{Proposition}
\label{ibunnot}
By passing to homology, the composite
\begin{eqnarray*}
\GGG_{\CCC}\to \GGG_{\CCC}[[t]] &\stackrel{\chi}{\to}& \End(HH_*(\CCC/k))[[t]]\oplus \End(HH_*(\CCC/k))((t))[-1] \\  &\to& (\End(HH_*(\CCC/k))((t))/\End(HH_*(\CCC/k))[[t]])[-1]
\end{eqnarray*}
induces $H^*(\frac{I_{\CCC}}{t}):H^*(\GG_{\CCC})\simeq H^*(\GGG_{\CCC})\to (\End(HH_*(\CCC/k))((t))/\End(HH_*(\CCC/k))[[t]])[-1]$.
\end{Proposition}

\Proof
Recall first that
the isomorphism $(\GGG_{\CCC}\otimes_kk[\eta])((t))\simeq \GGG_{\CCC}((t))\oplus \GGG_{\CCC}((t))[1]$
in Lemma~\ref{mysteriousisom}
is given by
$\GGG_{\CCC}((t))\oplus (\GGG_{\CCC}\otimes \eta) ((t))\to\GGG_{\CCC}((t))\oplus  (\GGG_{\CCC}\otimes \eta) ((t)) $ which sends $(g_1\cdot t^n, g_2\otimes \eta\cdot t^m)$
to $(g_1\cdot t^n,g_2\cdot t^{m+1}[1])$ (we use the notation in the proof of
Lemma~\ref{mysteriousisom}).
In particular, we identify $(\GGG_{\CCC}\otimes \eta) ((t))$ with
$\GGG_{\CCC} ((t))[1]$ by the assignment
$g\otimes\eta\cdot t^m \mapsto g\cdot t^{m+1}[1]$.
Note that $\widehat{\mathbb{A}}_{\CCC}^L:\GGG_{\CCC}\oplus (\GGG_{\CCC}\otimes\eta)\to  \End(HH_*(\CCC/k))$
carries $(g_1, g_2\otimes \eta)$
to 
 $\mathbb{A}_{\CCC}^L(g_1)+\mathbb{I}_{\CCC}(g_2)$.
Thus, 
\[
\GGG_{\CCC}((t))\oplus \GGG_{\CCC}((t))[1]\simeq\GGG_{\CCC}((t))\oplus (\GGG_{\CCC}\otimes\eta)((t)) \stackrel{ \widehat{\mathbb{A}}_{\CCC}^L((t))}{\longrightarrow}  \End(HH_*(\CCC/k))((t))
\]
carries $(g_1\cdot t^n, g_2\cdot t^m[1])$
to 
 $\mathbb{A}_{\CCC}^L(g_1)\cdot t^n+ \frac{\mathbb{I}_{\CCC}(g_2)}{t}\cdot t^{m}.$
Taking into account
the definition of $J$ in the proof of Lemma~\ref{extensiontokernel}, we see that
the composite $(\GGG_{\CCC}\otimes \eta) ((t))\simeq \GGG_{\CCC} ((t))[1]\stackrel{J\otimes_{k[[t]]}k((t))}{\to} \End(HH_*(\CCC/k))((t))$
is given by $\frac{\mathbb{I}_{\CCC}}{t}$.
By Proposition~\ref{desiredcompatibility}
and Corollary~\ref{kernel2}, $H^*(\mathbb{I}_{\CCC})=H^*(I_{\CCC})$.
Thus, by Proposition~\ref{desiredcompatibility} and Corollary~\ref{explicitformula},
$\GGG_{\CCC}\to  (\End(HH_*(\CCC/k))((t))/\End(HH_*(\CCC/k))[[t]])[-1]$
induces 
$H^*(\frac{I_{\CCC}}{t})=H^*(J):H^*(\GGG_{\CCC})\simeq H^*(\GG_{\CCC}) \to (\End(HH_*(\CCC/k))((t))/\End(HH_*(\CCC/k))[[t]])[-1]$.
\QED

\section{Proof}

\label{completion}

In this section we will complete the proof of the main result of this paper,
see Theorem~\ref{generalBTT}.
As a consequence we obtain Theorem~\ref{CYBTT} about Calabi-Yau
$\infty$-categories.

\subsection{}
\label{GBTTsection}
Let $\CCC$ be a small $k$-linear stable $\infty$-category.
As in Section~\ref{Periodmap}, in Section~\ref{GBTTsection}
we assume the $\SO(2)$-action of $\HH_\bullet(\CCC/k)$ is trivial.
We also assume that $\HH^\bullet(\CCC/k)$ and $\HH_\bullet(\CCC/k)$ are bounded
with finite dimensional cohomology,
i.e., $HH_n(\CCC/k)=0$ for $|n|>>0$ and $\dim_k HH_n(\CCC/k)<\infty$ for $n\in \ZZ$. 
When $\CCC$ is a
smooth and proper Calabi-Yau $\infty$-category, these conditions are satisfied.

\begin{Lemma}
\label{simpleendo}
There exists a canonical isomorphism
$\End(HH_\ast(\CCC/k))[[t]]\simeq  \End_{k[t]}(HH_*(\CCC/k)\otimes_kk[t])$.
Moreover, $\End(HH_\ast(\CCC/k))((t))\simeq  \End_{k[t^{\pm1}]}(HH_*(\CCC/k)\otimes_kk[t^{\pm1}])$.
\end{Lemma}

\Proof
In this proof, we use the condition that $HH_*(\CCC/k)$ is bounded.
By this assumption, the graded vector space
$\End(HH_\ast(\CCC/k))$ is bounded. 
Note that $\End(HH_*(\CCC/k))$ has bounded amplitude
so that 
\[
\End(HH_\ast(\CCC/k))[[t]]\simeq \End(HH_\ast(\CCC/k))\otimes_kk[t]\simeq \End_{k[t]}(HH_*(\CCC/k)\otimes_kk[t]).
\]
The latter isomorphism follows from the first isomorphism.
\QED

\begin{Theorem}
\label{generalBTT}
Suppose the following condition: the linear map
\[
\displaystyle HH^s(\CCC/k)\to \bigoplus_{i\in \ZZ}\Hom_{k}(HH_i(\CCC/k),HH_{i-s}(\CCC/k))
\]
given by $a\mapsto H^*(I_{\CCC})(a)$
is injective for any integer $s$.
Then $\GG_{\CCC}$ is quasi-abelian.
Namely,
 it is quasi-isomorphic to a dg Lie algebra with the zero bracket.
\end{Theorem}

\Proof
According to Proposition~\ref{ibunnot} and Proposition~\ref{mainreduction3},
it is enough to prove that $H^*(\frac{I_{\CCC}}{t})$ is injective.
By Lemma~\ref{simpleendo},
\[
\End(HH_*(\CCC/k))((t))\simeq \End_{k[t^{\pm1}]}(HH_*(\CCC/k)\otimes_kk[t^{\pm1}])\simeq \Hom_{k}(HH_*(\CCC/k),HH_*(\CCC/k)\otimes_k k[t^{\pm1}])
\]
where $\Hom_k(-,-)$ indicates the space of $k$-linear maps. Thus, there exists a canonical isomorphism
\[
\End_{k[t^{\pm1}]}(HH_*(\CCC/k))((t))/\End_{k[t]}(HH_*(\CCC/k))[[t]]\simeq 
\oplus_{i\in \ZZ,j\in \ZZ, r<0}\Hom_k(HH_i(\CCC/k), HH_j(\CCC/k)\cdot t^r)
\]
such that the (homological) degree of elements in $\Hom_k(HH_i(\CCC/k), HH_j(\CCC/k)\cdot t^r)$
is $-2r+j-i$.
It follows that our condition amounts to the injectivity of $H^*(\frac{I_{\CCC}}{t})$.
\QED

\subsection{}
We will prove that if $\CCC$ has a Calabi-Yau structure, then
the condition in Theorem~\ref{generalBTT} is satisfied (see Corollary~\ref{CYBTT}).

We first  recall the notion of 
left Calabi-Yau structures on a $k$-linear stable $\infty$-category $\CCC$
(see e.g. \cite{Gi}, \cite{BD} for details).
Let $\mathcal{HN}_{\bullet}(\CCC/k)$ be the negative cyclic complex
(computing negative cyclic homology) which is defined to be
the homotopy fixed points $\HH_\bullet(\CCC/k)^{\SO(2)}$.
Let $\mathcal{HN}_{\bullet}(\CCC/k)=\HH_\bullet(\CCC/k)^{\SO(2)}\to \HH_\bullet(\CCC/k)$ be the canonical morphism in $\Mod_k$.
Let $k[n]\to \mathcal{HN}_{\bullet}(\CCC/k)$
be a morphism in $\Mod_k$ that can be identified with an element $\xi \in \mathcal{HN}_{n}(\CCC/k)$ of the degree $n$.
Let $\xi_0$ be the image of $\xi$ in $HH_n(\CCC/k)=\HH_n(\CCC/k)$,
which we regard as a morphism $k[n]\to \HH_\bullet(\CCC/k)$.

We write $\CCC^e=\CCC^{op}\otimes_k\CCC$
for the tensor product in $\ST_k$.
We refer to an object of
$\Ind(\CCC^e)\simeq \Fun^{\textup{ex}}((\CCC^e)^{op}, \SP)$
as a right $\CCC^e$-module.
Here $\Fun^{\textup{ex}}((\CCC^e)^{op},\SP)$ is the full subcategory
of $\Fun((\CCC^e)^{op},\SP)$, which consists of exact functors.
By the Morita theory,
there exists a canonical equivalence $\Ind(\CCC^{op}\otimes_k\CCC)\simeq  \mathcal{M}or_k(\Ind(\CCC),\Ind(\CCC))$ where $\mathcal{M}or_k(\Ind(\CCC),\Ind(\CCC))$
is the internal hom/mapping object
in $\PR_k$ (cf. Section~\ref{Pr3}) whose underlying $\infty$-category classifies $k$-linear functors (i.e., morphisms of
$\Mod_k^\otimes$-modules).
We denote by $\CCC_m$ the (right) $\CCC^e$-module
corresponding to the identity functor $\Ind(\CCC)\to \Ind(\CCC)$.
The object $\CCC_m$ is usually referred to as
the diagonal right $\CCC^e$-module
because when
$\CCC\simeq \Mod_A$, $\CCC_m$ is $A$ having the diagonal
$A$-$A$-bimodule structure.
The permutation symmetry $\CCC^{op}\otimes_k\CCC\simeq \CCC\otimes_k\CCC^{op}$ determines the
equivalence $\Fun^{\textup{ex}}((\CCC^e)^{op},\SP)\simeq \Fun^{\textup{ex}}(\CCC^e,\SP)$.
The image of $\CCC_m$ in $\Fun^{\textup{ex}}(\CCC^e,\SP)$
is a left $\CCC^e$-module, which we also denote by $\CCC_m$.

Let $\Hom_{\CCC^e}^r(-,-)$ and $\Hom_{\CCC^e}^l(-,-)$ denote
the hom/mapping complex in $\Fun^{\textup{ex}}((\CCC^e)^{op}, \SP)$
and $\Fun^{\textup{ex}}(\CCC^e, \SP)$, respectively.
These are defined as objects of $\Mod_k$.
Let $\otimes_{\CCC^e}$
be the relative tensor product of the right $\CCC^e$-module
and the left $\CCC^e$-module.
For a left $\CCC^e$-module $\mathcal{M}$, there exists
the tautological morphism
$\textup{Hom}_{\CCC^e}^l(\mathcal{M},\mathcal{M})\otimes_k\mathcal{M}\to \mathcal{M}$ as left $\CCC^e$-modules.
This morphism determines
\[
\sigma_{\CCC,\mathcal{M}}:\textup{Hom}_{\CCC^e}^l(\mathcal{M},\mathcal{M})\otimes_k(\CCC_m\otimes_{\CCC^e}\mathcal{M})\to (\CCC_m\otimes_{\CCC^e}\mathcal{M}).
\]
The explicit definition of $\CCC_m\otimes_{\CCC^e}\mathcal{M}$
using spectral categories will be given below.
As mentioned above, $\CCC_m$
can also be considered as a left $\CCC^e$-module.
When $\mathcal{M}=\CCC_m$, 
there exist equivalences $\Hom^l_{\CCC^e}(\CCC_m,\CCC_m)\simeq \HH^\bullet(\CCC/k)$
and $\CCC_m\otimes_{\CCC^e}\CCC_m \simeq \HH_\bullet(\CCC/k)$
in $\Mod_k$.
Thus, $\sigma_{\CCC,\mathcal{M}}$
defines 
$\sigma_{\CCC}:\HH^\bullet(\CCC/k)\otimes \HH_\bullet(\CCC/k)\simeq \textup{Hom}^l_{\CCC^e}(\CCC_m,\CCC_m)\otimes_k(\CCC_m\otimes_{\CCC^e}\CCC_m) \stackrel{\sigma_{\CCC,\CCC_m}}{\to} \CCC_m\otimes_{\CCC^e}\CCC_m=\HH_\bullet(\CCC/k)$.

Suppose that $\CCC$ is smooth over $k$ (see Section~\ref{Pr3}).
For the left $\CCC^e$-module $\CCC_m$, we denote by $\CCC_m^!$
the inverse dualizing module defined as a right $\CCC^e$-module.
The inverse dualizing module $\CCC_m$
can be characterized as
a right $\CCC^e$-module
which corepresents
the functor $\Fun^{\textup{ex}}((\CCC^e)^{op},\SP)\to \SSS$ given by
$\mathcal{M}\mapsto \Omega^{\infty}(\mathcal{M}\otimes_{\CCC^e}\CCC_m)$
where $\SSS$ is the $\infty$-category of $\infty$-groupoids/spaces,
and $\Omega^{\infty}(\mathcal{M}\otimes_{\CCC^e}\CCC_m)$ indicates the underlying space through $\Mod_k\to \SP\to \SSS$ (cf. \cite[Example 2.11]{BD}).
There is a functorial equivalence 
$\Hom^r_{\CCC^e}(\CCC_m^!,\mathcal{M})\simeq \mathcal{M}\otimes_{\CCC^e}\CCC_m$.
In particular, $\Hom_{\CCC^e}^r(\CCC_m^!,\CCC_m)\simeq \HH_\bullet(\CCC/k)$.

\begin{Definition}
Suppose that $\CCC$ is smooth over $k$.
We say that $\xi$ is a left Calabi-Yau structure on $\CCC$ (or simply a Calabi-Yau structure
on $\CCC$)
if $\xi_0:k[n]\to \HH_\bullet(\CCC/k)\simeq \Hom^r_{\CCC^e}(\CCC_m^!,\CCC_m)$ determines an equivalence $\CCC^!_m\to \CCC_m[-n]$.
\end{Definition}

\begin{Remark}
There exists a canonical equivalence $\CCC^!_m\otimes_{\CCC^e}\CCC_m\simeq \textup{Hom}_{\CCC^e}^l(\CCC_m,\CCC_m)$ (see \cite[(2.11)]{BD}).
Moreover, it is straightforward to see that $\sigma_{\CCC}$
can be identified with
$(\CCC^!_m\otimes_{\CCC^e}\CCC_m)\otimes_k\Hom^r_{\CCC^e}(\CCC_m^!,\CCC_m) \to \CCC_m\otimes_{\CCC^e}\CCC_m$
determined by the tautological morphism $\Hom_{\CCC^e}^r(\CCC_m^!,\CCC_m)\otimes \CCC_m^!\to \CCC_m$.
Thus if $\xi_0$ defines an equivalence $\CCC_m^!\to \CCC_m[-n]$,
the composite
\[
d_{\CCC}:\HH^\bullet(\CCC/k)\otimes_kk[n] \stackrel{\textup{id}\otimes\xi_0}{\longrightarrow} \HH^\bullet(\CCC/k)\otimes \HH_\bullet(\CCC/k)\stackrel{\sigma_{\CCC}}{\longrightarrow}\HH_\bullet(\CCC/k)
\]
is an equivalence in $\Mod_k$
(i.e., which induces an isomorphism $HH^*(\CCC/k)\simeq HH_*(\CCC/k)[-n])$
of graded vector spaces).
\end{Remark}

\begin{Definition}
Let $D_{\CCC}$ be the composite
\[
\HH^\bullet(\CCC/k)\otimes_kk[n] \stackrel{\textup{id}\otimes\xi_0}{\longrightarrow} \HH^\bullet(\CCC/k)\otimes \HH_\bullet(\CCC/k)\stackrel{I_{\CCC}}{\longrightarrow}\HH_\bullet(\CCC/k).
\]
\end{Definition}

\begin{Proposition}
\label{dualitycoincide}
The morphism $D_{\CCC}:\HH^{\bullet}(\CCC/k)[n] \to \HH_\bullet(\CCC/k)$
is equivalent to the morphism $d_{\CCC}: \HH^{\bullet}(\CCC/k)[n]\to \HH_\bullet(\CCC/k)$.
In particular, $D_{\CCC}$ is an equivalence if and only if $d_{\CCC}$ is an equivalence.
\end{Proposition}

\begin{Corollary}
\label{CYBTT}
Let $\CCC$ be a $k$-linear proper Calabi-Yau stable $\infty$-category. 
Then $\GG_{\CCC}$ is quasi-abelian.
\end{Corollary}

\Proof
By the degeneration of the Hodge-de Rham spectral sequence for $\CCC$
(\cite{Kal0}, \cite{Ma}),
the $\SO(2)$-action on $\HH_\bullet(\CCC/k)$ is trivial
for a smooth and proper stable $\infty$-category $\CCC$.
By the Calabi-Yau structure, $d_\CCC$ is an equivalence.
It follows from Proposition~\ref{dualitycoincide} that
$D_{\CCC}$ is an equivalence.
Thus, the condition in Theorem~\ref{generalBTT} is satisfied.
Then Theorem~\ref{generalBTT} implies our assertion.
\QED

Proposition~\ref{dualitycoincide} was proved in \cite[Theorem 8.2]{ID}
in the case when $\CCC$ is equivalent to the stable $\infty$-category
of $A$-module spectra (or dg $A$-modules).
We will deduce Proposition~\ref{dualitycoincide} from \cite[Theorem 8.2]{ID}.
To this end, before the proof we review Hochschild homology
spectra in terms of symmetric spectra (see \cite{HSS} for symmetric spectra).
Let $\SPS$ be the symmetric monoidal category of symmetric spectra.
Let $\KKK$ be a commutative ring symmetric spectrum that represents the base field $k$.
Let $\SPS(\KKK)$ be the symmetric monoidal category of $\KKK$-module symmetric spectra.
It admits a stable symmetric
monoidal $\KKK$-model structure defined in \cite[Theorem 2.6]{conv}
whose weak equivalences are stable equivalences.

Let $\CAT_{\KKK}$ be the category of $\KKK$-spectral categories,
i.e., categories enriched over the monoidal category $\SPS(\KKK)$.
Let $\CAT_{\KKK}^{pc}$ be the full subcategory of $\CAT_{\KKK}$
spanned by those $\KKK$-spectral categories such that
every mapping spectrum
is cofibrant in $\SPS(\KKK)$ (we refer to such a category as
a pointwise-cofibrant $\KKK$-spectral category).
Let $\CAT_{\KKK}^{pc}[M^{-1}]$ be the $\infty$-category
obtained by inverting Morita equivalences (see e.g. \cite{BGT1} for Morita equivalences).
By \cite[Theorem 5.1]{Co} (see also \cite[Propoition 6.7]{ID}),
there exists a categorical equivalence $\CAT_{\KKK}^{pc}[M^{-1}]\simeq \ST_k$.
Using this equivalence we take a
$\KKK$-spectral category $\CC$ which represents $\CCC$.
We will assume that $\CC$ is pretriangulated (see \cite[4.4]{BM} for the notion of pretriangulated spectral categories).
By taking a replacement  \cite[Propoition 6.3]{ID},
we may and will assume that $\CC$ is pointwise-cofibrant.

Let $\Fun_{\SPS(\KKK)}(\CC^{op}\otimes\CC,\SPS(\KKK))$
denote the category of left $\CC^{op}\otimes \CC$-modules.
By a left $\CC^{op}\otimes \CC$-module we mean
a $\SPS(\KKK)$-enriched functor 
$\CC^{op}\otimes\CC \to \SPS(\KKK)$.
By \cite[A.1.1]{SS}, $\Fun_{\SPS(\KKK)}(\CC^{op}\otimes\CC,\SPS(\KKK))$
has a combinatorial $\SPS(\KKK)$-enriched model structure
in which weak equivalence are objectwise stable equivalences
and fibrations are objectwise fibrations in $\SPS(\KKK)$
(referred to as the projective model structure).
By inverting weak equivalences we obtain the localized $\infty$-category
$\Fun_{\SPS(\KKK)}(\CC^{op}\otimes\CC,\SPS(\KKK))[W^{-1}]$ 
which is equivalent to
$\Ind((\CCC^e)^{op})\simeq \Fun^{\textup{ex}}(\CCC^e,\SP)$.
Let $\mathbb{M}$ be a left $\CC^{op}\otimes \CC$-module which represents
$\mathcal{M}\in \Fun^{\textup{ex}}(\CCC^e,\SP)$.
Assume that $\mathbb{M}(X,X')$
are cofibrant in $\SPS(\KKK)$ for any $X,X'\in \CC$
(take a cofibrant left $\CC^{op}\otimes\CC$-module with respect to
the projective model structure).
The Hochschild-Mitchell simplicial object $\HH^{\Delta}_{\bullet}(\CC, \mathbb{M})$ is defined as a simplicial diagram in $\SPS(\KKK)$ of the form
\[
\xymatrix{
\cdots\displaystyle \bigoplus_{X_0,X_1,X_2} \mathbb{M}(X_2,X_0)\otimes \CC(X_0,X_1)\otimes \CC(X_1,X_2)\ar@<1.0ex>[r] \ar[r] \ar@<-1.0ex>[r] & 
\displaystyle \bigoplus_{X_0,X_1} \mathbb{M}(X_1,X_0)\otimes \CC(X_0,X_1) \ar@<0.5ex>[r] \ar@<-0.5ex>[r] \ar@<-0.5ex>[l] \ar@<0.5ex>[l]& \displaystyle \bigoplus_{X_0\in \CC} \mathbb{M}(X_0,X_0) \ar[l]
}
\]
where $\otimes$ indicates the smash product over $\KKK$, and $\oplus$ indicates
the wedge sum.
The $p$-th term is defined to be
\[
\displaystyle \bigoplus_{X_0,\ldots,X_p} \mathbb{M}(X_p,X_0)\otimes \CC(X_0,X_1)\otimes\cdots \otimes \CC(X_{p-1},X_p).
\]
The degeneracy maps are induced by the composition in $\CC$
and the left $\CC^{op}\otimes\CC$-module structure in $\mathbb{M}$
in the standard way.
The face map are induced by the unit maps $\KKK\to \CC(X_i,X_i)$.
The geometric realization computes $\CCC_m\otimes_{\CCC^e}\mathcal{M}$.
We refer to \cite[Section 3.1, Section 5, Section 6]{BM} for the details of Hochschild-Mitchell simplicial objects including the Morita invariance
and the invariance under weak equivalences of $\mathbb{M}$.
If $\mathbb{M}$ is $\CC_m:\CC^{op}\otimes\CC\to \SPS(\KKK)$
given by $(X,Y)\mapsto \CC(X,Y)$,
the geometric realization (i.e., the homotopy colimit) 
of $\HH_\bullet^{\Delta}(\CC):=\HH_\bullet^{\Delta}(\CC,\CC_m)$ is the Hochschild homology
spectrum.

Let $\ZZZ$ be a cofibrant object in $\SPS(\KKK)$ which represents 
$\HH^\bullet(\CCC/k)$.
By taking a replacement we suppose that
$\CC_m$ is a cofibrant and fibrant left $\CC^{op}\otimes \CC$-module
which represents $\CCC_m$.
Consider the tautological morphism
$t:\HH^\bullet(\CCC/k)\otimes \CCC_m\to \CCC_m$
of left $\CCC^{op}\otimes\CCC$-modules (i.e., a morphism in $\Fun^{\textup{ex}}(\CCC^{op}\otimes_k\CCC,\SP)$).
Let $\ZZZ\otimes \CC_m$ denote the tensor with $\ZZZ$.
Since the category of left $\CC^{op}\otimes \CC$-modules
is a $\SPS(\KKK)$-enriched model category, $\ZZZ\otimes \CC_m$ is cofibrant.
Let $t':\ZZZ\otimes \CC_m\to \CC_m$ be a
morphism of left $\CC^{op}\otimes\CC$-modules
which represents $t$.
Then $t'$ induces a morphism $\sigma^{\Delta}_{\CC}:\ZZZ\otimes \HH_\bullet^{\Delta}(\CC)\to \HH_\bullet^{\Delta}(\CC)$ of simplicial diagrams, whose $p$-term is
given by the wedge sum of
\[
(\ZZZ\otimes\CC(X_p,X_0))\otimes \CC(X_0,X_1)\otimes\cdots \otimes \CC(X_{p-1},X_p)\stackrel{t'\otimes\textup{id}}{\to} \CC(X_p,X_0)\otimes \CC(X_0,X_1)\otimes\cdots \otimes \CC(X_{p-1},X_p).
\]
The geometric realization of $\sigma_{\CC}^{\Delta}:\ZZZ\otimes \HH_\bullet^{\Delta}(\CCC/k)\to \HH_\bullet^{\Delta}(\CCC/k)$
rerpesents $\sigma_{\CCC}$.

{\it Proof of Proposition~\ref{dualitycoincide}.}
Let $\{\CCC_{\alpha}\}_{\alpha\in \Lambda}$
be the poset consisting of those idempotent-complete
stable subcategories $\CCC_{\alpha}$ in $\CCC$ such that
$\Ind(\CCC_{\alpha})$ admits a single compact generator
(so that $\CCC_{\alpha}\simeq \RPerf_{A_{\alpha}}$ for some dg algebra $A_{\alpha}$).
Since $\HH_\bullet(-/k)$ commutes with filtered colimits (cf. \cite[10.2]{BGT1}),
the equivalence $\colim_{\alpha\in \Lambda}\CCC_{\alpha}\simeq \CCC$
induces $\colim_{\alpha\in \Lambda}\HH_\bullet(\CCC_{\alpha}/k)\simeq \HH_\bullet(\CCC/k)$.
The morphism $\xi_0:k[n]\to \HH_\bullet(\CCC/k)$
factors as $\xi_0:k[n]\to \HH_\bullet(\CCC_{\alpha}/k)\to  \HH_\bullet(\CCC/k)$ for some $\alpha$.
Let $\xi_0'\in \HH_n(\CCC_{\alpha}/k)$ be the associated morphism
(note that $\xi_0$ is the image of $\xi_0'$ under $HH_n(\CCC_{\alpha}/k)\to HH_n(\CCC/k)$).
Note that $I_{\CCC}$ is restricted/promoted to 
$I_{\CCC,\CCC_{\alpha}}:\HH^\bullet(\CCC/k)\otimes \HH_\bullet(\CCC_{\alpha}/k) \to \HH_\bullet(\CCC_{\alpha}/k)$
(see Remark~\ref{restrictedaction}) which makes the diagram
\[
\xymatrix{
\HH^\bullet(\CCC/k)\otimes \HH_\bullet(\CCC_{\alpha}/k) \ar[r]^(0.6){I_{\CCC,\CCC_{\alpha}}} \ar[d] & \HH_\bullet(\CCC_{\alpha}/k)  \ar[d] \\
\HH^\bullet(\CCC/k)\otimes \HH_\bullet(\CCC/k) \ar[r]^(0.6){I_{\CCC}} & \HH_\bullet(\CCC/k)  
}
\]
commute
where the vertical arrows is induced by $\HH_\bullet(\CCC_{\alpha}/k)\to \HH_\bullet(\CCC/k)$.
It will suffice to show that $\HH^\bullet(\CCC/k)\otimes_kk[n] \stackrel{\textup{id}\otimes\xi'_0}{\longrightarrow} \HH^\bullet(\CCC/k)\otimes \HH_\bullet(\CCC_{\alpha}/k)\stackrel{I_{\CCC,\CCC_{\alpha}}}{\longrightarrow} \HH_\bullet(\CCC_{\alpha}/k)\to \HH_\bullet(\CCC/k)$ is equivalent to
$d_{\CCC}$.
Let $\CC$ be the pretriangulated pointwise-cofibrant spectral category over $\KKK$
which represents $\CCC$ (defined as above).
Let $\CC_{\alpha}$ be the spectral category
that has a single object representing a compact generator of $\Ind(\CCC_{\alpha})$.
Then we have the left commutative diagram
\[
\xymatrix{
\ZZZ\otimes \HH_\bullet^\Delta(\CC_{\alpha}) \ar[r]^{\sigma^{\Delta}_{\CC,\CC_{\alpha}}} \ar[d] & \HH_\bullet^{\Delta}(\CC_{\alpha})  \ar[d] & & \HH^\bullet(\CCC/k)\otimes \HH_\bullet(\CCC_{\alpha}/k) \ar[r]^(0.6){\sigma_{\CCC,\CCC_\alpha}} \ar[d] & \HH_\bullet(\CCC_{\alpha}/k)  \ar[d]\\
\ZZZ\otimes \HH_\bullet^\Delta(\CC) \ar[r]^{\sigma_{\CC}^\Delta} & \HH_\bullet^\Delta(\CC)   & & \HH^\bullet(\CCC/k)\otimes \HH_\bullet(\CCC/k) \ar[r]^(0.6){\sigma_{\CCC}} & \HH_\bullet(\CCC/k) 
}
\]
where $\sigma^{\Delta}_{\CC,\CC_{\alpha}}$ is defined by restricting
$\sigma_{\CC}^{\Delta}$ to $\CC_{\alpha}$.
The right square
diagram is obtained by taking geometric realizations of the left square
diagram of simplicial objects.
The comparison result \cite[Theorem 8.2]{ID} (and the proof)
says that $I_{\CCC,\CCC_{\alpha}}$
is equivalent to $\sigma_{\CCC,\CCC_\alpha}$.
It follows that $\HH^\bullet(\CCC/k)\otimes \HH_\bullet(\CCC_{\alpha}/k)\stackrel{I_{\CCC,\CCC_{\alpha}}}{\longrightarrow} \HH_\bullet(\CCC_{\alpha}/k)\to \HH_\bullet(\CCC/k)$ is equivalent to the composite
$\HH^\bullet(\CCC/k)\otimes \HH_\bullet(\CCC_{\alpha}/k)\to \HH^\bullet(\CCC/k)\otimes \HH_\bullet(\CCC/k) \stackrel{\sigma_{\CCC}}{\to} \HH_\bullet(\CCC/k)$.
Thus our assertion follows.
\QED


\end{document}